\documentclass[11pt,a4paper]{article}
\usepackage[T1]{fontenc}
\usepackage{lmodern,microtype,amsmath,amssymb,amsthm,mathtools,mathrsfs}
\usepackage[margin=27mm]{geometry}
\usepackage{enumitem}
\usepackage[hidelinks]{hyperref}
\hypersetup{%
  pdftitle={The maximum product of nontrivial cross t-intersecting families},
  pdfauthor={Biao Wu and Huajun Zhang}%
}
\newtheorem{theorem}{Theorem}[section]
\newtheorem{lemma}[theorem]{Lemma}
\newtheorem{proposition}[theorem]{Proposition}

\newtheorem{problem}{Problem}
\theoremstyle{remark}
\newcommand{\calH}{\mathcal H}\newcommand{\calG}{\mathcal G}
\numberwithin{equation}{section}
\setlist{topsep=3pt,itemsep=2pt,parsep=0pt,leftmargin=*}
\title{The maximum product of nontrivial\\cross $t$-intersecting families}
\author{Biao Wu\thanks{School of Mathematics and Statistics, Hunan Normal
University, Changsha 410081, P. R. China.
Email: \href{mailto:wu@hunnu.edu.cn}{\texttt{wu@hunnu.edu.cn}}.}
\qquad Huajun Zhang\thanks{Department of Mathematics, Shaoxing University,
Shaoxing 312000, P. R. China. Corresponding author.
Email: \href{mailto:huajunzhang@zjnu.cn}{\texttt{huajunzhang@zjnu.cn}}.}}
\date{}
\begin{document}
\maketitle
\begin{abstract}
Let $\mathcal A$ and $\mathcal B$ be nonempty families of $k$-subsets
of $\{1,\ldots,n\}$ such that $|A\cap B|\ge t$ for all $A\in\mathcal A$
and $B\in\mathcal B$, and suppose that each family has fewer than $t$
common elements. For
$k\ge t+2\ge4$ and $n\ge(t+1)(k-t+1)$, we prove
\[
 |\mathcal A||\mathcal B|\le\max\{h_t(n,k)^2,b_t(n,k)^2\},
\]
where
\[
 h_t(n,k)=\binom{n-t}{k-t}-\binom{n-k-1}{k-t}+t,\quad
 b_t(n,k)=(t+2)\binom{n-t-2}{k-t-1}+\binom{n-t-2}{k-t-2}.
\]
The bound is sharp, and we characterize all extremal pairs, whose
families have equal sizes but need not coincide. The proof combines
clique lifts, exchanges of generating families and a density inequality
for cross-intersecting uniform layers. This gives a product version
of the Hilton--Milner--Frankl theorem throughout the Erd\H{o}s--Ko--Rado
range. In particular, it answers a question of Wen and Lv on whether
a threshold of order $kt$ suffices for the sharp product bound.
\end{abstract}

{\small\noindent\textbf{Keywords:} cross $t$-intersecting families;
nontrivial families; generating families; Hilton--Milner theorem.\par}

\section{Introduction}
Let $[n]=\{1,\ldots,n\}$, and let $\binom{[n]}k$ denote the family
of all $k$-subsets of $[n]$. A family $\mathcal F\subseteq\binom{[n]}k$
is \emph{$t$-intersecting} if $|F\cap F'|\ge t$ for all
$F,F'\in\mathcal F$. Two families $\mathcal A,\mathcal B\subseteq
\binom{[n]}k$ are \emph{cross $t$-intersecting} if $|A\cap B|\ge t$
for all $A\in\mathcal A$ and $B\in\mathcal B$.
For $t=1$, we say \emph{intersecting} and \emph{cross-intersecting},
respectively. For integers $u\ge0$ and $v$, we take $\binom uv=0$
unless $0\le v\le u$.

The Erd\H{o}s--Ko--Rado theorem determines the maximum size of a
$t$-intersecting family for sufficiently large $n$. Its sharp form is
as follows.

\begin{theorem}[Erd\H{o}s--Ko--Rado~\cite{EKR}, Frankl~\cite{FranklEKR},
Wilson~\cite{Wilson}]\label{EKR}
Let $1\le t\le k\le n$ and $n\ge(t+1)(k-t+1)$.
If $\mathcal F\subseteq\binom{[n]}k$ is $t$-intersecting, then
\[
 |\mathcal F|\le\binom{n-t}{k-t}.
\]
\end{theorem}

The bound is attained by a \emph{$t$-star}, the family of all $k$-sets
containing a fixed $t$-set. Erd\H{o}s, Ko and Rado proved the bound
for sufficiently large $n$, including the sharp range $n\ge2k$ when
$t=1$. Frankl established the sharp threshold for $t\ge15$, and
Wilson proved it for all $t$. For $t<k$, the threshold
$(t+1)(k-t+1)$ cannot be lowered. Below this threshold, requiring
many points from a small fixed block can outperform a $t$-star.
Ahlswede and Khachatrian~\cite{AKCIT} determined the maximum for
all parameters in their complete intersection theorem. Thus the
Erd\H{o}s--Ko--Rado range is precisely the natural range in which
a full $t$-star provides the benchmark for the extremal problem.

To exclude families contained in a $t$-star, a nonempty family
$\mathcal F$ is called \emph{nontrivial} if
$|\bigcap_{F\in\mathcal F}F|<t$.
The resulting problem is to determine the largest nontrivial
$t$-intersecting family.
Hilton and Milner~\cite{HM} answered this question for $t=1$.

\begin{theorem}[Hilton--Milner~\cite{HM}]\label{thm:classical-hm}
Let $k\ge3$ and $n>2k$. If $\mathcal F\subseteq\binom{[n]}k$ is a
nonempty intersecting family with $\bigcap_{F\in\mathcal F}F=\varnothing$,
then
\[
 |\mathcal F|\le
 \binom{n-1}{k-1}-\binom{n-k-1}{k-1}+1.
\]
\end{theorem}

This theorem identifies the first obstruction to star structure:
a family larger than the displayed bound must have a common point.
Frankl~\cite{F78-2} determined the corresponding maximum for
nontrivial $t$-intersecting families when $n$ is sufficiently large.
Ahlswede and Khachatrian~\cite{AK} subsequently determined the
sharp bound in the Erd\H{o}s--Ko--Rado range, where two constructions
determine the maximum. One modifies a $t$-star by
adding a small set of exceptional members; the other requires
$t+1$ points from a fixed $(t+2)$-set.

Let $T,X\subseteq[n]$ be disjoint, with $|T|=t$ and $|X|=k-t+1$,
and let $Z\subseteq[n]$ have size $t+2$. Define
\begin{align*}
 \mathcal H(T,X)
 &=\{F\in\tbinom{[n]}k:T\subseteq F,\ F\cap X\ne\varnothing\}
   \cup\{(T\setminus\{p\})\cup X:p\in T\},\\
 \mathcal A_1(Z)
 &=\{F\in\tbinom{[n]}k:|F\cap Z|\ge t+1\}.
\end{align*}
These are nontrivial $t$-intersecting families of sizes
\[
 h_t(n,k)=\binom{n-t}{k-t}-\binom{n-k-1}{k-t}+t,\quad
 b_t(n,k)=(t+2)\binom{n-t-2}{k-t-1}+\binom{n-t-2}{k-t-2},
\]
respectively. We refer to $\mathcal H(T,X)$ as the Hilton--Milner
construction and to $\mathcal A_1(Z)$ as the radius-one construction.
Put $m_t(n,k)=\max\{h_t(n,k),b_t(n,k)\}$.

\begin{theorem}[Hilton--Milner--Frankl~\cite{HM,F78-2,AK}]
\label{thm:nontrivial-ak}
Let $k\ge t+2\ge4$ and $n\ge(t+1)(k-t+1)$.
If $\mathcal F\subseteq\binom{[n]}k$ is nontrivial and
$t$-intersecting, then
\[
 |\mathcal F|\le m_t(n,k).
\]
Equality holds if and only if $\mathcal F$ is a family
$\mathcal H(T,X)$ or $\mathcal A_1(Z)$ of size $m_t(n,k)$.
\end{theorem}

For $n>(t+1)(k-t+1)$, the equality statement follows from the complete
nontrivial-intersection theorem~\cite{AK}. At the boundary, the
complete intersection theorem~\cite{AKCIT} permits only a $t$-star
or a radius-one family, and nontriviality excludes the star.

Cross $t$-intersection imposes the intersection condition only on
pairs of sets drawn from opposite families, so the bound for a
single $t$-intersecting family need not apply to either side. The sum and product of the family sizes can have
different extremal structures; see~\cite{XiangZhang}. For the product,
two copies of an extremal one-family construction give a natural
lower bound. The question is whether unequal families can give a
larger product.

Pyber~\cite{Pyber} initiated the product problem for cross-intersecting
families and proved the star-product bound
$|\mathcal A||\mathcal B|\le\binom{n-1}{k-1}^2$ for $n\ge2k$.
Matsumoto and Tokushige~\cite{MatsumotoTokushige} extended this
bound to different uniformities. For general $t$,
Tokushige~\cite{TokushigeEigen} used eigenvalue methods to prove
$|\mathcal A||\mathcal B|\le\binom{n-t}{k-t}^2$ under
$k/n<1-2^{-1/t}$ and conjectured that the Erd\H{o}s--Ko--Rado
threshold suffices. Frankl, Lee, Siggers and Tokushige~\cite{FLST}
proved the bound for $t\ge14$ and $n\ge(t+1)k$.
The authors~\cite{ZhangWuCrossT} reached the sharp threshold for
$t\ge3$, and Tanaka and Tokushige~\cite{TanakaTokushige} settled
$t=2$ by semidefinite programming. Together, these results give
the following product analogue of Theorem~\ref{EKR}.

\begin{theorem}[Pyber~\cite{Pyber}, the authors~\cite{ZhangWuCrossT},
Tanaka--Tokushige~\cite{TanakaTokushige}]\label{thm:cross-ekr}
Let $k\ge t+1\ge2$ and $n\ge(t+1)(k-t+1)$.
If $\mathcal A,\mathcal B\subseteq\binom{[n]}k$ are cross
$t$-intersecting, then
\[
 |\mathcal A||\mathcal B|\le\binom{n-t}{k-t}^2.
\]
\end{theorem}

Two full $t$-stars with the same centre attain this bound. It is
therefore natural to ask for the maximum product after excluding
star containment on each side, just as the Hilton--Milner theorem
does for a single family.

We require each family to be nontrivial:
\begin{equation}\label{eq:separate-nontriviality}
 \left|\bigcap_{A\in\mathcal A}A\right|<t,
 \qquad
 \left|\bigcap_{B\in\mathcal B}B\right|<t.
\end{equation}
We call such a pair \emph{nontrivial}. The weaker condition
$|\bigcap_{F\in\mathcal A\cup\mathcal B}F|<t$ still allows one
family to lie in a $t$-star. He, Li and the authors~\cite{HeLiWuZhang}
studied that problem for $t\ge3$ and different uniformities, and
Bao, Ji and Yu~\cite{BaoJiYu} settled $t=2$.
Condition~\eqref{eq:separate-nontriviality} excludes these one-sided
star configurations, giving the product version of the
Hilton--Milner--Frankl problem.

For $t=1$, the natural candidate is the square of the Hilton--Milner
number,
\[
 h_1(n,k)=\binom{n-1}{k-1}-\binom{n-k-1}{k-1}+1.
\]
Frankl and Wang~\cite{FranklWangHM} proved the sharp product bound
$h_1(n,k)^2$ for $k\ge8$ and $n\ge4k$.
Frankl~\cite{FranklThree} settled $k=3$, and Frankl and
Wang~\cite{FranklWangII} extended their result to $n\ge2k+1$
for $k\ge8$. Independently, Huang~\cite{Huang} obtained the bound for all
$k\ge3$ and $n\ge2k$. The authors~\cite{WZ} determined the maximum
product for different uniformities in the range
$n\ge2\ell\ge2k\ge6$ and $n>k+\ell$.

For $t\ge2$, Frankl and Wang~\cite{FranklWang} proved the bound
$m_t(n,k)^2$ for $n\ge4(t+2)^2k^2$ and $k\ge5$.
Wen and Lv~\cite[Theorem~1.7]{WL} improved the range to
$n\ge t+\max\{15(t+2)^2(k-t),10k^2\}$ for $k\ge t+2\ge4$.
They showed that the maximum is $h_t(n,k)^2$ when $k\ge2t+2$
and $b_t(n,k)^2$ when $k\le2t+1$, and determined the extremal pairs.
In particular, the sharp product bound in their range is
$\max\{h_t(n,k)^2,b_t(n,k)^2\}$. Wen and Lv asked whether a condition
of order $kt$ suffices for this bound. We state the question in our
notation.

\begin{problem}[{Wen and Lv~\cite[Problem~1]{WL}}]\label{prob:wen-lv}
Does there exist an absolute constant $C>0$ such that, whenever
$k\ge t+2\ge4$ and $n\ge Ckt$, any two nonempty, nontrivial cross
$t$-intersecting families $\mathcal A,\mathcal B\subseteq\binom{[n]}k$
satisfy
\[
 |\mathcal A||\mathcal B|\le\max\{h_t(n,k)^2,b_t(n,k)^2\}?
\]
\end{problem}

We determine the maximum product throughout
$n\ge(t+1)(k-t+1)$. The maximum is the larger of the two construction
values, and the extremal families need not be equal. To state the
latter possibility, let $T,X,Y\subseteq[n]$ satisfy
$|T|=t$, $|X|=|Y|=k-t+1$ and $T\cap(X\cup Y)=\varnothing$, and define
\[
 \mathcal H(T;X,Y)=
 \{F\in\tbinom{[n]}k:T\subseteq F,\ F\cap Y\ne\varnothing\}
 \cup\{(T\setminus\{p\})\cup X:p\in T\}.
\]
Then $\mathcal H(T;X,X)=\mathcal H(T,X)$ and
$|\mathcal H(T;X,Y)|=h_t(n,k)$.
The pair $(\mathcal H(T;X,Y),\mathcal H(T;Y,X))$ is cross
$t$-intersecting if and only if $|X\cap Y|\ge2$: exceptional members
on opposite sides omitting different points of $T$ have intersection
size $t-2+|X\cap Y|$, and the other pairs of members satisfy the
required bound when $|X\cap Y|\ge2$. In this case we call it a
\emph{dual Hilton--Milner pair}, as in~\cite{WL}.

\begin{theorem}\label{thm:main}
Let $k\ge t+2\ge4$ and $n\ge(t+1)(k-t+1)$.
If $\mathcal A,\mathcal B\subseteq\binom{[n]}k$ are nonempty,
nontrivial cross $t$-intersecting families, then
\[
 |\mathcal A||\mathcal B|\le m_t(n,k)^2.
\]
Equality holds if and only if, up to a simultaneous permutation of
$[n]$, one of the following holds:
\begin{itemize}
 \item[\rm(i)] $b_t(n,k)\ge h_t(n,k)$ and
 $\mathcal A=\mathcal B=\mathcal A_1(Z)$ for a $(t+2)$-set $Z$;
 \item[\rm(ii)] $h_t(n,k)\ge b_t(n,k)$ and
 $\mathcal A=\mathcal H(T;X,Y)$, $\mathcal B=\mathcal H(T;Y,X)$,
 where $|T|=t$, $|X|=|Y|=k-t+1$, $T\cap(X\cup Y)=\varnothing$
 and $|X\cap Y|\ge2$.
\end{itemize}
\end{theorem}

When $h_t(n,k)=b_t(n,k)$, both alternatives occur. Every extremal
pair has equal family sizes, although its two families may differ.
At $n=(t+1)(k-t+1)$, we have
\[
 b_t(n,k)=\binom{n-t}{k-t},\qquad
 b_t(n,k)-h_t(n,k)=\binom{t(k-t)}{k-t}-t>0,
\]
so only \textup{(i)} occurs. Theorem~\ref{thm:main} answers
Problem~\ref{prob:wen-lv} affirmatively with $C=3/2$, since
\[
 (t+1)(k-t+1)\le(t+1)k\le\tfrac32kt
 \qquad(t\ge2).
\]

Theorem~\ref{thm:main} gives a product version of the
Hilton--Milner--Frankl theorem throughout the Erd\H{o}s--Ko--Rado range.
The maximum is the square of the one-family maximum, with dual
Hilton--Milner pairs giving the additional equality cases.

Simultaneous shifting preserves cardinalities and cross
$t$-intersection, but it may give one family a common $t$-set.
We therefore apply shifts only while both families remain nontrivial.
The process ends at a shifted pair or at an unshifted terminal pair
for which every changing shift makes at least one family trivial.
\emph{Critical shiftedness}, defined for a single family in
Section~\ref{sec:critical-family}, turns this obstruction to further
compression into a structural restriction. The authors~\cite{WZ} exploit
a decomposition into parts with disjoint cores to study nontrivial
cross-intersecting products. The second author~\cite{ZhangNontrivialSum} uses
local generator reductions in the critically shifted case for
nontrivial cross-intersecting sums. An analogous stopping condition
under a covering-number constraint is used by
the second author~\cite{ZhangCoverThree} to study large maximal intersecting
families with covering number three.

We treat both terminal cases using clique lifts, whose weights
record the higher clique layers that recover a cross $t$-intersecting
family from an ordinary cross-intersecting family. In the critical
case, deleting a common $(t-1)$-set and applying the two-core
decomposition of the authors~\cite{WZ} reduces the problem to a
weighted density estimate. Its strictness excludes unshifted
terminal pairs from equality. In the shifted case, standard
generating families~\cite{XiangZhang} and paired exchanges extending
those in~\cite{ZhangWuCrossT} preserve nontriviality and reduce the
problem to kernels controlled by the same density estimate.
At equality, these reductions force a common $t$-transversal of
size $t+1$; reversing the shifts recovers the original extremal pairs.

Section~\ref{sec:clique-reduction} develops the clique lifts,
generator exchanges and common-transversal estimates. The critical
case is proved in Section~\ref{sec:critical-family}. Section~\ref{sec:shifted}
treats the shifted case, including the necessary kernel estimates.
Section~\ref{sec:equality} determines the equality cases and completes
the proof of Theorem~\ref{thm:main}. The algebraic estimates for the
exchanges are collected in Appendix~\ref{sec:exchange-numerics}.

\section{Preliminaries and structural reductions}
\label{sec:clique-reduction}

\subsection{Clique lifts and their companions}
Put $a=k-t+1$ and $N=n-t+1$. For an $a$-uniform family
$\mathcal H$ on an $N$-element set $V$, let
\[
 \mathcal K_q(\mathcal H)
 =\{Y\in\tbinom Vq:\tbinom Ya\subseteq\mathcal H\},
 \qquad c_q(\mathcal H)=|\mathcal K_q(\mathcal H)|
 \quad(q\ge a).
\]
Thus $\mathcal K_q(\mathcal H)$ consists of the $q$-element cliques of
$\mathcal H$. For a $(t-1)$-set $C$ disjoint from $V$, define the
\emph{clique lift} and its weight by
\begin{align*}
 L_C(\mathcal H)
 &=\{(C\setminus D)\cup Y:D\subseteq C,
                         \ Y\in\mathcal K_{a+|D|}(\mathcal H)\},\\
 W_t(\mathcal H)
 &=\sum_{d=0}^{t-1}\binom{t-1}{d}c_{a+d}(\mathcal H).
\end{align*}
We identify $C\mathbin{\dot\cup}V$ with $[n]$.
Then $L_C(\mathcal H)$ is $k$-uniform and has size $W_t(\mathcal H)$.
The members containing $C$ form the base layer
$C\vee\mathcal H=\{C\cup H:H\in\mathcal H\}$.

For $\mathcal A\subseteq\binom{[n]}k$, write
\[
 \mathcal A^{\perp_t}
 =\{B\in\tbinom{[n]}k:|A\cap B|\ge t\text{ for every }A\in\mathcal A\}
\]
for the family of all $k$-sets that $t$-intersect every member of
$\mathcal A$. Similarly, $\mathcal H^\perp$ denotes the family of all
members of $\binom Va$ meeting every member of $\mathcal H$. We call
these the full companions of $\mathcal A$ and $\mathcal H$, respectively.

\begin{lemma}\label{lem:clique-lift-duality}
If $\mathcal H,\mathcal G\subseteq\binom Va$ are cross-intersecting,
then $L_C(\mathcal H)$ and $L_C(\mathcal G)$ are cross $t$-intersecting.
Moreover,
\begin{equation}\label{eq:lift-duality}
 (L_C(\mathcal H))^{\perp_t}=L_C(\mathcal H^\perp).
\end{equation}
If $\mathcal H$ is nonempty, then $L_C(\mathcal H)$ is nontrivial if
and only if $\bigcap\mathcal H=\varnothing$.
\end{lemma}
\begin{proof}
For $Y\in\binom V{a+e}$, we have
\begin{equation}\label{eq:lift-companion-test}
 Y\in\mathcal K_{a+e}(\mathcal H^\perp)
 \quad\Longleftrightarrow\quad
 |Y\cap H|\ge e+1\quad\text{for every }H\in\mathcal H.
\end{equation}
Indeed, an $a$-subset of $Y$ avoids $H$ exactly when
$|Y\cap H|\le e$.

If $X\in\mathcal K_{a+d}(\mathcal H)$ and
$Y\in\mathcal K_{a+e}(\mathcal G)$, then
\begin{equation}\label{eq:cross-clique-intersection}
 |X\cap Y|\ge d+e+1.
\end{equation}
Otherwise partition $X\cap Y$ into parts of sizes at most $d$ and $e$,
delete them from $X$ and $Y$, respectively, and complete the deletions
to sizes $d$ and $e$. This gives disjoint members of $\mathcal H$ and
$\mathcal G$. For $D,E\subseteq C$ with $|D|=d$ and $|E|=e$, it follows that
\[
 |((C\setminus D)\cup X)\cap((C\setminus E)\cup Y)|
 \ge t-1-|D\cup E|+d+e+1\ge t.
\]
This proves the first assertion and gives
$L_C(\mathcal H^\perp)\subseteq(L_C(\mathcal H))^{\perp_t}$.
Conversely, let $(C\setminus E)\cup Y$ be a member of
$(L_C(\mathcal H))^{\perp_t}$. Its intersection with every base member
$C\cup H$ has size at least $t$, so $|Y\cap H|\ge |E|+1$ for every
$H\in\mathcal H$. The reverse inclusion now follows from
\eqref{eq:lift-companion-test}.

If $\bigcap\mathcal H=\varnothing$, the base layer has intersection
$C$, so the lift is nontrivial. If every $H\in\mathcal H$ contains
$x$, there are no cliques of order greater than $a$: every larger set
has an $a$-subset omitting $x$. The lift is then $C\vee\mathcal H$,
and every member contains the $t$-set $C\cup\{x\}$.
\end{proof}

\subsection{Shifting and saturation}
A cross $t$-intersecting pair is \emph{saturated} if neither family
can be enlarged while preserving cross $t$-intersection; equivalently,
the families are full companions of one another.
A pair of nonempty nontrivial families of maximum product is saturated,
since adding members cannot destroy nontriviality.

For $i<j$, the shift $s_{ij}$ replaces a member $A$ by $A-j+i$
when $j\in A$, $i\notin A$, and the replacement is absent from the
family. A family fixed by every such shift is \emph{shifted}.
A pair is \emph{terminal} if every simultaneous left shift that changes
it gives one of the two families a common $t$-set.

\begin{theorem}\label{thm:simplified-weighted-reduction}
Let $t\ge2$, $k\ge t+2$ and $n>2k-t$. The maximum product of the
sizes of two nonempty nontrivial cross $t$-intersecting families in
$\binom{[n]}k$ is
\[
 \max W_t(\mathcal H)W_t(\mathcal G),
\]
where the maximum is over nonempty cross-intersecting
$\mathcal H,\mathcal G\subseteq\binom{[N]}a$ with
$\bigcap\mathcal H=\bigcap\mathcal G=\varnothing$.
Every maximizing pair can be transformed by simultaneous shifts,
preserving both cardinalities, into $(L_C(\mathcal H),L_C(\mathcal G))$ with
$\mathcal G=\mathcal H^\perp$ and $\mathcal H=\mathcal G^\perp$.
\end{theorem}
\begin{proof}
By Lemma~\ref{lem:clique-lift-duality}, the lifts of any pair counted
on the right are nontrivial cross $t$-intersecting families. It remains
to show that every pair of maximum product can be brought to this form.

A simultaneous shift preserves both cardinalities and cross
$t$-intersection. Indeed, if both sets of a pair move or neither
moves, their intersection size is unchanged. If $A$ moves and $B$
does not, the intersection can decrease only when $j\in B$ and
$i\notin B$. In this case $B-j+i$ was in the original family, and its
intersection with the original $A$ has the same size as that of $B$
with the shifted $A$.

Starting from a maximizing pair, apply changing shifts that preserve
nontriviality. Every intermediate pair still maximizes the product
and is therefore saturated. The sum of the element labels strictly
decreases, so the process ends at a terminal pair
$(\mathcal A,\mathcal B)$.

\emph{Case 1: both families are shifted.}
Take $C=[t-1]$ and $V=[n]\setminus C$. Let
$\mathcal H=\{H\in\binom Va:C\cup H\in\mathcal A\}$ and
$\mathcal G=\{G\in\binom Va:C\cup G\in\mathcal B\}$.
If $(C\setminus D)\cup Y\in\mathcal A$, then for every
$Z\in\binom Ya$ we may shift the points of $Y\setminus Z$ into $D$.
Thus $C\cup Z\in\mathcal A$, and
\[
 \mathcal A\subseteq L_C(\mathcal H),\qquad
 \mathcal B\subseteq L_C(\mathcal G).
\]
The inclusions show that $\mathcal H$ and $\mathcal G$ are nonempty.
They are cross-intersecting, since $C\cup H\in\mathcal A$ and
$C\cup G\in\mathcal B$ share only $t-1$ points in $C$.
By Lemma~\ref{lem:clique-lift-duality}, their lifts are cross
$t$-intersecting. Saturation thus gives equality in both
inclusions. The same lemma shows that $\bigcap\mathcal H=\varnothing$
and $\bigcap\mathcal G=\varnothing$.

\emph{Case 2: a shift gives one family a common $t$-set.}
Suppose $s_{ij}(\mathcal A)$ has a common $t$-set. A shift can add
only $i$ to the common intersection. If $i$ is already common it
changes nothing; if $j$ is common but $i$ is not, it replaces $j$
by $i$ without increasing the intersection size. Hence
$C=\bigcap\mathcal A$ has size $t-1$ and avoids $i,j$.
Delete $C$ to obtain an ordinary nontrivial family $\mathcal H$ whose
image under $s_{ij}$ has common point $i$. Every member of $\mathcal H$
meets $\{i,j\}$, and $\mathcal H$ cannot contain both $P+i$ and $P+j$
for an $(a-1)$-set $P$ avoiding $i,j$, since $P+j$ would survive the
shift. An $(a+1)$-clique would either contain such a pair or have an
$a$-subset avoiding both $i$ and $j$. Hence there are no cliques of
size greater than $a$, and
\[
 \mathcal A=C\vee\mathcal H=L_C(\mathcal H),\qquad
 \mathcal B=L_C(\mathcal H^\perp)
\]
by saturation and \eqref{eq:lift-duality}. Set
$\mathcal G=\mathcal H^\perp$. Since its lift is nonempty and
nontrivial, $\mathcal G$ is nonempty and has empty intersection.

In either case, saturation and \eqref{eq:lift-duality}, followed
by comparison of the base layers, give
$\mathcal G=\mathcal H^\perp$ and $\mathcal H=\mathcal G^\perp$.
\end{proof}

\subsection{Generating families and terminal exchanges}\label{sec:generators}
Following Xiang and the second author~\cite{XiangZhang}, we use the standard
generating family. For a nonempty $\mathcal F\subseteq\binom{[n]}k$, put
\[
 \mathscr P(\mathcal F)=\{E\subseteq[n]:|E|\le k,
       \{K\in\tbinom{[n]}k:E\subseteq K\}\subseteq\mathcal F\},
 \qquad
 \mathscr G(\mathcal F)=\min_{\subseteq}\mathscr P(\mathcal F).
\]
Thus $\mathscr G(\mathcal F)$ is the antichain of minimal sets whose
complete stars lie in $\mathcal F$. We also call its members the
canonical generators of $\mathcal F$. A $k$-set belongs to
$\mathcal F$ if and only if it contains a member of
$\mathscr G(\mathcal F)$.

\begin{lemma}\label{lem:canonical-prefix}
Suppose that $\mathcal F$ is shifted. If $E\in\mathscr P(\mathcal F)$,
$j\in E$, $i<j$ and $i\notin E$, then
$E-j+i\in\mathscr P(\mathcal F)$ and hence $E-j+i$ contains a member
of $\mathscr G(\mathcal F)$.
For a nonempty generator $E$, define
\[
 D(E)=\{K\in\tbinom{[n]}k:K\cap[\max E]=E\}.
\]
If $\varnothing\notin\mathscr G(\mathcal F)$, then
\[
 \mathcal F=\mathop{\dot\bigcup}_{E\in\mathscr G(\mathcal F)}D(E),
 \qquad |D(E)|=\binom{n-\max E}{k-|E|}.
\]
\end{lemma}
\begin{proof}
Take a $k$-set $K$ containing $E-j+i$. If $j\in K$, then
$E\subseteq K$ and $K\in\mathcal F$. If $j\notin K$, then
$K-i+j\in\mathcal F$ because it contains $E$. Since $\mathcal F$
is shifted, $K\in\mathcal F$ in this case as well.

For $K\in\mathcal F$, choose a generator $E\subseteq K$ with
$\max E$ as small as possible. If $x\in K\cap[\max E]\setminus E$,
then $x<\max E$, and the first assertion gives a generator contained
in $E-\max E+x$. Its maximum is smaller, a contradiction. Thus
$K\in D(E)$. If $D(E)\cap D(F)\ne\varnothing$ and
$\max E\le\max F$, then $E\subseteq F$; the antichain property
gives $E=F$. The count follows by choosing the points above $\max E$.
\end{proof}

\begin{lemma}
\label{lem:canonical-pairing}
Let $n>2k-t$, and let $\mathcal A,\mathcal B\subseteq\binom{[n]}k$
be nonempty shifted families that form a saturated cross
$t$-intersecting pair. Then $\mathscr G(\mathcal A)$ and
$\mathscr G(\mathcal B)$ are cross $t$-intersecting.
With
\[
 s=\max\{\max E:E\in\mathscr G(\mathcal A)\cup
                                  \mathscr G(\mathcal B)\},
\]
every $E\in\mathscr G(\mathcal A)$ with $\max E=s$ has a partner
$F\in\mathscr G(\mathcal B)$ satisfying
\[
 s\in E\cap F,\qquad |E\cap F|=t,\qquad
 E\cup F=[s],\qquad |E|+|F|=s+t.
\]
The same statement holds after interchanging the families. In
particular both generator systems have the same maximal support $s$,
and $s\le2k-t$.
If both families are nontrivial, then every generator has size at
least $t+1$.
\end{lemma}
\begin{proof}
Take $E\in\mathscr G(\mathcal A)$ and $F\in\mathscr G(\mathcal B)$.
Extend them to $k$-sets by first using points outside their union,
without increasing their intersection. If the union becomes $[n]$
before both sets reach size $k$, complete each set using points in
the other. The resulting intersection has size
$\max\{|E\cap F|,2k-n\}$. Since $n>2k-t$, cross $t$-intersection
forces $|E\cap F|\ge t$.

By minimality, $E-s$ is not a complete-star generator of $\mathcal A$.
Choose a $k$-set $K\notin\mathcal A$ containing $E-s$. Since the pair
is saturated, some $B\in\mathcal B$ satisfies $|K\cap B|<t$. Choose
$F\in\mathscr G(\mathcal B)$ contained in $B$. Then $|(E-s)\cap F|<t$.
Since $|E\cap F|\ge t$, we have $s\in F$ and $|E\cap F|=t$.
If $x\in[s]\setminus(E\cup F)$, Lemma~\ref{lem:canonical-prefix} gives
an element of $\mathscr G(\mathcal A)$ contained in $E-s+x$ and
meeting $F$ in at most $t-1$ points, a contradiction. Hence
$E\cup F=[s]$, and the remaining assertions follow by counting.

The generators have size at least $t$ because their families are
cross $t$-intersecting. If $E\in\mathscr G(\mathcal A)$ has size $t$,
then $E$ is contained in every generator of $\mathcal B$ and hence
in every member of $\mathcal B$, contrary to nontriviality. The same
argument applies to generators of $\mathcal B$.
\end{proof}

\begin{lemma}\label{lem:terminal-deletion}
Let $1\le t\le k<n$, and let $\mathcal F$ be a shifted nontrivial
family. Suppose that its canonical generators are supported on $[s]$
and have size at least $t$. Delete all $j$-element generators
containing $s$, and let $\mathcal F^-$ be the family generated by
those remaining. Then $\mathcal F^-$ is nonempty, and the members of
$\mathcal F\setminus\mathcal F^-$ are precisely those in the prefix
cylinders of the deleted generators. If all members of $\mathcal F^-$
contain a common $t$-set, then $j=s-1$.
\end{lemma}
\begin{proof}
If no generator is deleted, there is nothing to prove. A deleted
generator $E$ cannot be $[s]$: it would then be the only generator,
contrary to nontriviality. Shifting $s$ to a point of $[s-1]\setminus E$
gives a set containing a generator supported on $[s-1]$. This generator
is retained, so $\mathcal F^-$ is nonempty.
If a member containing a deleted generator $E$ has an additional
prefix point $x\in[s-1]\setminus E$, it also contains a retained
generator inside $E-s+x$. Conversely, a member of $D(E)$ cannot
contain another generator, because every generator is supported on
$[s]$ and the generators form an antichain. Thus precisely the
cylinders $D(E)$ of the deleted generators disappear.

Let $C$ be a common $t$-set of $\mathcal F^-$. Since $k<n$, the
intersection of all $k$-sets containing a given generator is that
generator. Applying this to the retained generator in $[s-1]$ gives
$C\subseteq[s-1]$.
Since $\mathcal F$ is nontrivial, some deleted generator $E$ omits a
point $c\in C$. For every $x\in[s-1]\setminus E$, shifting $s$ to $x$
gives a surviving generator and hence $C\subseteq E-s+x$.
It follows that every such $x$ equals $c$. Therefore
$E=[s]\setminus\{c\}$ and $j=s-1$.
\end{proof}

\begin{lemma}\label{lem:retained-canonical}
Suppose $\mathscr G(\mathcal B)=\mathscr G_0\dot\cup\mathscr E$,
where every member of $\mathscr E$ has size $k$, and let $\mathcal B_0$
be generated by $\mathscr G_0$. Then
$\mathscr G(\mathcal B_0)=\mathscr G_0$.
\end{lemma}
\begin{proof}
Every retained generator remains minimal, since
$\mathcal B_0\subseteq\mathcal B$. Conversely, a complete-star
generator $H$ of $\mathcal B_0$ contains an original canonical generator
$E$ of $\mathcal B$. If $E\in\mathscr E$, then $|E|=k$ forces $H=E$,
but $E\notin\mathcal B_0$ by the original antichain property. Hence
$E\in\mathscr G_0$, proving the assertion.
\end{proof}

The following form of the exchange in
\cite[Lemma 2.3]{ZhangWuCrossT} records the size changes and the
conditions preserving nontriviality.

\begin{lemma}\label{lem:paired-terminal-exchange}
Let $n>2k-t$, and let $\mathcal A$ and $\mathcal B$ be nonempty
shifted families that form a saturated cross $t$-intersecting pair.
Suppose that their canonical generators have common support $[s]$.
Fix $i,j$ with $i+j=s+t$ and put
\[
\mathcal I=\{E\in\mathscr G(\mathcal A):s\in E,\ |E|=i\},\qquad
\mathcal J=\{F\in\mathscr G(\mathcal B):s\in F,\ |F|=j\}.
\]
Suppose $\alpha=|\mathcal I|>0$ and $\beta=|\mathcal J|>0$.
Let $\mathcal A^+$ be generated by $\mathscr G(\mathcal A)$ together
with $\{E\setminus\{s\}:E\in\mathcal I\}$, and let $\mathcal A^-$
be generated by $\mathscr G(\mathcal A)\setminus\mathcal I$.
Define $\mathcal B^+,\mathcal B^-$ analogously using $\mathcal J$.
Both $(\mathcal A^+,\mathcal B^-)$ and
$(\mathcal A^-,\mathcal B^+)$ are cross $t$-intersecting, and
\begin{align*}
|\mathcal A^+|&=|\mathcal A|+\alpha\binom{n-s}{k-i+1},&
|\mathcal A^-|&=|\mathcal A|-\alpha\binom{n-s}{k-i},\\
|\mathcal B^+|&=|\mathcal B|+\beta\binom{n-s}{k-j+1},&
|\mathcal B^-|&=|\mathcal B|-\beta\binom{n-s}{k-j}.
\end{align*}
Moreover,
\begin{equation}\label{eq:terminal-shadow-bound}
|\mathcal A|\ge
\alpha\binom{n-s}{k-i}
+\alpha\frac{s-i}{i}\binom{n-s+1}{k-i},
\end{equation}
and the analogous inequality holds for $\mathcal B$.
If both original families are nontrivial and $i,j\ne s-1$, then
both families in each exchanged pair are nontrivial.
\end{lemma}
\begin{proof}
Suppose $E\in\mathcal I$ and $F\in\mathscr G(\mathcal B)$ satisfy
$|(E\setminus\{s\})\cap F|<t$.
Since $|E\cap F|\ge t$, we have $s\in E\cap F$ and $|E\cap F|=t$.
If $x\in[s]\setminus(E\cup F)$, shifting $s$ to $x$ in $E$
gives a complete-star generator meeting $F$ in fewer than $t$ points.
Hence $E\cup F=[s]$ and $|F|=s+t-i=j$, so $F\in\mathcal J$.
Deleting $\mathcal J$ therefore makes $\mathcal A^+$ and
$\mathcal B^-$ cross $t$-intersecting. The reverse exchange is symmetric.

Deletion removes exactly the terminal cylinders, by the argument
of Lemma~\ref{lem:terminal-deletion}. A new member containing
$E-s$ has trace exactly $E-s$ on $[s]$: a further prefix point
$x$ would give the original generator $E-s+x$. Conversely, such a
trace cannot contain an original generator, by minimality of $E$.
Thus addition creates disjoint cylinders of size
$\binom{n-s}{k-i+1}$, proving the four formulas.

Let $\mathcal I'=\{E\setminus\{s\}:E\in\mathcal I\}$ and let
$\nabla\mathcal I'$ be its rank-$i$ upper shadow in $[s-1]$.
Each member of $\mathcal I'$ has $s-i$ extensions and each shadow
member contains at most $i$ members of $\mathcal I'$, so
$|\nabla\mathcal I'|\ge\alpha(s-i)/i$.
Since $\mathcal A$ is shifted, every member of the shadow is a
complete-star generator. The $\alpha$ terminal cylinders have traces
of size $i-1$ on $[s-1]$, whereas the cylinders with traces in
$\nabla\mathcal I'$ have traces of size $i$. Counting these disjoint
collections gives~\eqref{eq:terminal-shadow-bound}.
Finally, enlargement preserves nontriviality, and both deletions
preserve it by Lemma~\ref{lem:terminal-deletion} when $i,j\ne s-1$.
\end{proof}

\subsection{Common transversals}
\label{subsec:common-transversal}
We call a set a common $t$-transversal if it meets every member of
both families in at least $t$ points. When this set has size $t+1$,
we use ordinary cross-intersection between the sections obtained by
fixing the trace on it. We first consider a section contained in a star.

\begin{lemma}\label{lem:offset-star-companion}
Let $q\ge3$, $r\ge2$, $N\ge qr$, and put
\[
 a=\binom N{r-1},\qquad b=\binom{N-1}{r-1},\qquad
 D=\binom Nr-\binom{N-r}r.
\]
If a nonempty family in $\binom{[N]}r$ of size $s$ is contained in a star,
and its full ordinary cross-intersecting companion has size $g$, then
\begin{equation}\label{eq:offset-star}
 (a+qs)(a+qg)\le\bigl(a+\max\{qb,D\}\bigr)^2.
\end{equation}
Equality holds if and only if $s=b$ and $qb\ge D$.
\end{lemma}
\begin{proof}
We first establish a density inequality. For cross-intersecting
$\alpha$- and $\beta$-uniform families on $[m]$, where $\alpha,\beta>0$ and
$m>\alpha+\beta$, let $u,v$ be their densities. Then
\begin{equation}\label{eq:transversal-density}
 u+v+\frac{uv}{L}\le1,
 \qquad L=\frac{\alpha\beta}{m(m-\alpha-\beta)}.
\end{equation}
By symmetry, assume $\alpha\le\beta$. Equip both uniform layers with
the uniform probability measure. For a function on $\beta$-sets,
let $T$ give its average over the $\beta$-sets disjoint from a given
$\alpha$-set.
Write $(x)_j=x(x-1)\cdots(x-j+1)$, with $(x)_0=1$.
For $f_J^{(b)}(E)=\mathbf1_{J\subseteq E}$ and $j=|J|\le\alpha$,
\[
 Tf_J^{(\beta)}=\frac{(\beta)_j}{(m-\alpha)_j}
       \sum_{I\subseteq J}(-1)^{|I|}f_I^{(\alpha)}.
\]
The adjoint identity interchanges $\alpha,\beta$, since choosing a
uniform $\alpha$-set and then a uniform disjoint $\beta$-set gives the
uniform distribution on disjoint pairs. Let $U_j$ be the span of
$\{f_J^{(\alpha)}:|J|\le j\}$. The displayed formula and its adjoint
show that $TT^*$ preserves $U_j$ and acts on $U_j/U_{j-1}$ by the
scalar
$\frac{(\alpha)_j(\beta)_j}{(m-\alpha)_j(m-\beta)_j}$.
Since $TT^*$ is self-adjoint, it acts by this scalar on
$U_j\cap U_{j-1}^{\perp}$. These subspaces together with the constants
span the space of functions on the $\alpha$-sets. The eigenvalues
decrease with $j$, as $m>\alpha+\beta$. Thus the operator norm of $T$
on functions of mean zero is at most
$\rho=\sqrt{\alpha\beta/((m-\alpha)(m-\beta))}$.
Write $f,g$ for the indicator functions of the two families. By
cross-intersection, $\langle f,Tg\rangle=0$. Hence
$\langle f-u,T(g-v)\rangle=-uv$, and Cauchy--Schwarz gives
$uv\le\rho\sqrt{u(1-u)v(1-v)}$. Squaring, cancelling when $uv>0$,
and using $\rho^2=L/(1+L)$ proves~\eqref{eq:transversal-density};
the case $uv=0$ is immediate.

Now suppose that the given family is contained in the star at $1$.
Write $g=b+h$. The full companion contains every member of this star,
and $0\le h\le D-b$ since every companion member meets a fixed member
of the original family. The sets obtained by deleting $1$ from the
original members are cross-intersecting with the companion members
that omit $1$. Applying~\eqref{eq:transversal-density} to these two
families gives
\[
 \frac hb\le\frac{C(1-U)}{L+U},\qquad
 U=\frac sb,\quad L=\frac{r(r-1)}{(N-1)(N-2r)},\quad
 C=\frac{N-r}{r}L.
\]
Set $\alpha=a/b=N/(N-r+1)$, $d=D/b$ and $z=N/r$.
Then $0<U\le1$ and $1<d<z$. The normalized product is at most
$(\alpha+qU)\left(\alpha+q+ q\min\left\{d-1,\frac{C(1-U)}{L+U}\right\}\right)$.
The branches meet at $U_0=L(z-d)/(L(z-1)+d-1)\in(0,1)$.
The first branch is increasing; the second, denoted $F(U)$, satisfies
\[
 F''(U)=\frac{2qC(1+L)(\alpha-qL)}{(L+U)^3}>0,
 \qquad qL<\frac1{q-2}\le1<\alpha.
\]
It suffices to bound $F(1)=(\alpha+q)^2$ and
$F(U_0)=(\alpha+qU_0)(\alpha+qd)$ by
$(\alpha+\max\{q,d\})^2$.

Suppose first that $r\ge3$. In the expansion
$D=\sum_{i=1}^r\binom{N-i}{r-1}$, the first term is $b$ and the
next two exceed $2b/3$ and $4b/9$. Thus $d>1+2/3+4/9>2$.
Also $\alpha\le q/(q-1)$ and
$L<1/(z(z-2))$, whence
\begin{equation}\label{eq:offset-crossing}
 U_0<\frac{z-d}{z-1+z(z-2)(d-1)}
       \le\frac1{4(d-1)^2+1}.
\end{equation}
For the last step put $w=z-d$: the denominator is
$w+(d-1)(w+d-1)^2$, and $(w+d-1)^2\ge4w(d-1)$.
If $2\le d\le q$, it suffices to show that
\[
 U_0\le\frac{q+\alpha(2-d)}{\alpha+qd}.
\]
The right side decreases with $\alpha$, so it is at least
$(q+1-d)/(1+(q-1)d)$. This lower bound is greater than the bound
in~\eqref{eq:offset-crossing}, since
\begin{align*}
 &(q+1-d)(4(d-1)^2+1)-(1+(q-1)d)\\
 &\qquad=(d-1)(4(d-1)(q+1-d)-q)>0.
\end{align*}
Here $(d-1)(q+1-d)\ge q-1$, by concavity on $[2,q]$.
If $d\ge q$, it suffices to show that
$U_0\le d(d-\alpha(q-2))/(q(\alpha+qd))$.
The right side decreases with $\alpha$ and is at least
\[
 \frac{d((q-1)d-q(q-2))}{q^2(1+(q-1)d)}
 \ge\frac{d}{q(1+(q-1)d)}\ge\frac1{q^2-q+1},
\]
whereas~\eqref{eq:offset-crossing} is strictly smaller than the last
quantity for $d\ge q\ge3$. The comparison at $U_0$ is strict.
Together with the monotonicity of the first branch and strict
convexity of the second, this proves~\eqref{eq:offset-star} for
$r\ge3$, with equality only at $U=1$ and $q\ge d$.

For $r=2$, $a=N$, $b=N-1$, $D=2N-3<qb$.
The companion sizes are $g=2N-3$ for $s=1$, $g=N$ for $s=2$,
and $g=N-1$ for $3\le s\le N-1$.
At $s=1,2$, the gaps from $(N+q(N-1))^2$ are respectively
$q^2(N-2)^2$ and $Nq(N-4)+q^2(N^2-4N+1)$, both positive
since $N\ge2q\ge6$. For $s\ge3$ the product increases strictly
to the same square at $s=N-1=b$. This also proves the equality
statement.
\end{proof}

\begin{theorem}\label{thm:transversal}
Let $t\ge2$, $k\ge t+2$ and $n\ge(t+1)(k-t+1)$. Suppose that
$\mathcal F,\mathcal G\subseteq\binom{[n]}k$ are nonempty nontrivial
cross $t$-intersecting families with a common $t$-transversal of
size $t+1$. Then
\[
 |\mathcal F||\mathcal G|\le\max\{h_t(n,k)^2,b_t(n,k)^2\}.
\]
\end{theorem}
\begin{proof}
Put $q=t+1$, $r=k-t$, $N=n-q\ge qr$, and retain $a,b,D$ from
Lemma~\ref{lem:offset-star-companion}. We have
$h_t(n,k)=a+D+q-1$ and $b_t(n,k)=a+qb$; write
$M=a+\max\{qb,D+q-1\}$.
Relabel the transversal as $T=[q]$ and add all $k$-sets containing $T$
to both families. The families remain cross $t$-intersecting and
nontrivial. In each family, the members containing $T$ contribute
$a$ sets. For each $i\in[q]$, remove $T\setminus\{i\}$ from the
members with this trace on $T$ to obtain a section of $r$-sets outside
$T$. Each family has at least two nonempty sections, since otherwise
it would have a common $t$-set in $T$. Sections from opposite families
with different omitted points must be cross-intersecting. Sections
with the same omitted point impose no further condition.

By Hilton's lexicographic lemma~\cite{HiltonLex}, we may replace
all sections by lexicographic initial segments of the same sizes on
$[n]\setminus T$. The required cross-intersections are preserved.
Let $\mathcal L(s)$ be the first $s$ outside $r$-sets in lexicographic
order, and let $g(s)$ be the size of its full companion. The full
companion is itself an initial segment: by Hilton's lemma, its
replacement is contained in the companion and has the same size.
Thus $g$ is nonincreasing, $g(b)=b$, and $g(1)=D$. Since
$\mathcal L(s)$ is contained in the full companion of its full
companion, $g(g(s))\ge s$ whenever $g(s)>0$.
Order the first family's section sizes so that the largest two are
$x\ge y\ge1$. The opposite sections have sizes at most $g(y)$
for the first label and $g(x)$ for every other label. Since at least
two opposite sections are nonempty, one of them has a different
label from the largest section. Hence $g(x)>0$ and $x\le D$.
Therefore
\begin{equation}\label{eq:transversal-product}
 |\mathcal F||\mathcal G|\le
 P(x,y):=(a+x+(q-1)y)(a+g(y)+(q-1)g(x)),
 \quad 1\le y\le x\le D.
\end{equation}

\emph{Both section sizes are at most $b$.} Suppose $y\le x\le b$.
We claim
$P(x,y)\le\max_{s\in\{x,y\}}(a+qs)(a+qg(s))$.
To verify this, put $d=x-y$, $e=g(y)-g(x)$,
$U=a+qx$, $V=a+qg(x)$. The three products are
$(U-(q-1)d)(V+e)$, $UV$, and $(U-qd)(V+qe)$.
The latter two agree at $e_0=dV/(U-qd)$.
At $e_0$ the first is at most $UV$, since
$(q-2)U-(q-1)^2d=(q-2)a-x+(q-1)^2y\ge0$;
here $x\le b\le a$ and $q\ge3$.
Below $e_0$ its difference from $UV$ increases with $e$;
above $e_0$ its difference from the other endpoint decreases, with
slope $x-(q-1)a-(q^2-q+1)y<0$.
This proves the claim (the case $d=0$ is immediate).
Each $\mathcal L(s)$ with $s\le b$ lies in a star, so
Lemma~\ref{lem:offset-star-companion} gives $P(x,y)\le M^2$.

\emph{Both section sizes are at least $b$.} Put $x'=g(y)$ and
$y'=g(x)$.
Then $1\le y'\le x'\le b$ and $g(g(s))\ge s$ gives
$P(x,y)\le P(x',y')$, reducing to the preceding case.

\emph{One section size is at least $b$ and the other is at most $b$.}
Suppose $x\ge b\ge y$, and put $X=x$ and $Y=g(y)$.
Both are at least $b$, while $g(X),g(Y)\le b$ and $g(Y)\ge y$.
Consequently
\begin{align*}
 P(x,y)&\le(a+X+(q-1)g(Y))(a+Y+(q-1)g(X))\\
 &\le(a+X+(q-1)g(X))(a+Y+(q-1)g(Y)).
\end{align*}
The second comparison follows by subtraction: the first product
minus the second is $(q-1)(X-Y)(g(X)-g(Y))\le0$.
For $z=X,Y$, form a family from all $k$-sets containing $T$, one
section $\mathcal L(z)$ and $q-1$ sections $\mathcal L(g(z))$.
The resulting family is $t$-intersecting, since
$\mathcal L(z)$ and $\mathcal L(g(z))$ are cross-intersecting, and
the sections $\mathcal L(g(z))$ are contained in a common outside
star. The intersection of the full $T$-star is $T$, and all omitted
labels have nonempty sections. Hence the family is nontrivial.
By Theorem~\ref{thm:nontrivial-ak}, $a+z+(q-1)g(z)\le M$.
Thus~\eqref{eq:transversal-product} is at most $M^2$ in every case.
\end{proof}

\section{Critically shifted families}\label{sec:critical-family}

Throughout this section put $a=k-t+1\ge3$ and $N=n-t+1$,
where $t\ge2$, and use the clique layers and weights from
Section~\ref{sec:clique-reduction}. Set
\[
 h=\binom{N-1}{a-1}-\binom{N-a-1}{a-1}+t,\quad
 b=(t+2)\binom{N-3}{a-2}+\binom{N-3}{a-3},\qquad B=\max\{h,b\}.
\]
Thus $h=h_t(n,k)$, $b=b_t(n,k)$ and $B=m_t(n,k)$.

In this section, a nonempty family on $[N]$ is \emph{nontrivial} if
the intersection of its members is empty. Following~\cite{WZ}, such
a family is \emph{critically shifted} if at least one left shift changes it and
every left shift that changes it makes it trivial.

Requiring every changing shift to make the family trivial yields
the two-core structure used below; one such shift alone does not
suffice. For a terminal pair, the family made trivial may depend on
the shift. The proof of Theorem~\ref{thm:critical} reduces this
pair condition to critical shiftedness of one family.

\begin{theorem}\label{thm:critical}
Let $a\ge3$, $t\ge2$, and
\[
 N\ge(t+1)(a-1)+2.
\]
Suppose that $\calH\subseteq\binom{[N]}a$ is critically shifted and that
$\calG\subseteq\binom{[N]}a$ is cross-intersecting with $\calH$. Then
\[
 |\calH|W_t(\calG)<B^2.
\]
Consequently, suppose that $\mathcal A,\mathcal B\subseteq\binom{[n]}k$
form a saturated terminal pair of nontrivial cross $t$-intersecting
families and that some simultaneous shift changes the pair. Then
\[
 |\mathcal A||\mathcal B|<\max\{h_t(n,k)^2,b_t(n,k)^2\}
 \qquad\text{for }n\ge(t+1)(k-t+1).
\]
\end{theorem}

By \cite[Lemma 3.4]{WZ}, every critically shifted $a$-uniform family
can be partitioned into two nonempty parts, one whose members contain
$R$ and the other whose members contain $S$, where $R\cap S=\varnothing$
and
\begin{equation}
 2\le |R|,|S|\le a,\qquad |R|+|S|\ge a+2.                 \label{eq:critical-core-structure}
\end{equation}
We also use the Erd\H{o}s--Ko--Rado product bound for cross
$j$-intersecting families of different uniformities. If
$u\ge v\ge j\ge1$, $m\ge(j+1)(u-j+1)$,
and $\mathcal U\subseteq\binom{[m]}u$,
$\mathcal V\subseteq\binom{[m]}v$ are cross $j$-intersecting, then
\begin{equation}
 |\mathcal U||\mathcal V|
 \le\binom{m-j}{u-j}\binom{m-j}{v-j}.                   \label{eq:unequal-rank-ekr}
\end{equation}
For $j=1$ this follows from Matsumoto and Tokushige~\cite{MatsumotoTokushige},
since $m\ge2u$. For $v>j$ and $j=2$, we use the result of
Bao, Ji and Xiang, as stated in Bao, Ji and Yu~\cite[Theorem~1.5]{BaoJiYu}.
For $v>j$ and $j\ge3$, we use He, Li and the
authors~\cite[Corollary~1.8]{HeLiWuZhang}.
The case $v=j$ is elementary. We may assume both families are nonempty.
Put $S=\bigcup\mathcal V$ and $s=|S|$. Every member of $\mathcal U$
contains $S$, so $j\le s\le u$ and
\[
 |\mathcal U||\mathcal V|\le\binom sj\binom{m-s}{u-s}.
\]
As $s$ increases from $j$ to $u$, the expression on the right is
nonincreasing: its consecutive ratio is
\[
 \frac{s+1}{s+1-j}\frac{u-s}{m-s}\le1,
\]
since $j(u-s)\le(s+1-j)(m-u)$ follows from
$m-u\ge j(u-j)+1$. Its value at $s=j$ proves
\eqref{eq:unequal-rank-ekr}.

We distinguish two cases. If both cores in
\eqref{eq:critical-core-structure} have size at least three, we bound
the product by counting the cliques according to their intersections
with the cores. Otherwise one core has size two and the other part
consists of a single set:
\begin{equation}
 \calH=(R\vee\mathcal P)\cup\{X\},\qquad
 |R|=2,\quad |X|=a,\quad R\cap X=\varnothing.             \label{eq:partial-star-form}
\end{equation}
Here $R\vee\mathcal P=\{R\cup P:P\in\mathcal P\}$.
We first bound the product when $\mathcal P$ contains all
$(a-2)$-sets outside $R$. We then use a density inequality to bound
the extra cliques that can occur when $\mathcal P$ is smaller.

\subsection{The full two-point star}\label{sec:endpoint}
\begin{lemma}\label{lem:endpoint}
Let $R,X\subseteq[N]$ be disjoint sets with $|R|=2$ and $|X|=a$. Put
\[
 \calH_* =\{H\in\tbinom{[N]}a:R\subseteq H\}\cup\{X\},
 \qquad\calG_* =\calH_*^\perp.
\]
If $N\ge(t+1)(a-1)+2$, then $|\calH_*|W_t(\calG_*)<h^2$.
\end{lemma}
\begin{proof}
Write $m=N-2$, $A=a-1$, and
$D=\binom m{a-2},\quad E=\binom m{a-1}-\binom{m-a}{a-1},\quad D_0=\binom{m-a}{a-2}$.
Since $N\ge2a$, every $a$-set disjoint from $R$ is disjoint from
some member of the full $R$-star. Thus an $a$-set belongs to
$\calG_*$ if and only if it meets both $R$ and $X$. It follows that
an $(a+d)$-set is a clique of $\calG_*$ if and only if it meets each
of $R$ and $X$ in at least $d+1$ points. In particular there are no
such cliques for $d\ge2$. Counting for $d=0,1$ gives
\[
 |\calH_*|=D+1,\quad
 Q:=W_t(\calG_*)=(t+1)E+D-[1+(t-1)a]D_0,\quad
 h=E+D-D_0+t.
\]
For $0\le j\le a$ put $q_j=\binom{m-j}{a-2}/D$, and let
$y=q_a$, $P=\sum_{j=1}^{A}q_j$, $u=1+P$.
Telescoping yields $E=D(P+y)$, whence
\begin{equation}
 h=Du+t,\qquad Q=Dw,\qquad
 w=(t+1)u-t-(t-1)Ay+y.                                \label{eq:formula-f1}
\end{equation}
The positive sequence $q_j$ is log-concave, since
$q_{j+1}/q_j=1-(a-2)/(m-j)$ decreases with $j$.
Consequently $q_iq_j\ge q_1q_{i+j-1}$ when $i+j-1\le A$.
Writing $Z=\sum_{r=1}^A rq_r$, we obtain
\[
 P^2\ge q_1Z,\qquad
 P-Ay=(a-2)\sum_{r=1}^A\frac{rq_r}{m-r}
       \le\frac{A-1}{m-A}Z.
\]
The quantity $c(m)=q_1-(t-1)(A-1)/(m-A)$ increases with $m$.
At $m=(t+1)A$ it is
$c(m)=\frac1A+\frac{A-1}{t(t+1)A}>\frac1A$.
Since $Z\ge A(A+1)y/2$ and $A\ge2$, it follows that
$P^2-(t-1)(P-Ay)\ge c(m)Z>\frac{A+1}{2}y>y$.
By \eqref{eq:formula-f1}, this is exactly $u^2-w>0$.
Also $w<(t+1)u\le2tu$, so
$h^2-(D+1)Q=D^2(u^2-w)+D(2tu-w)+t^2>0$.
\end{proof}

\subsection{A density inequality}\label{sec:density}
\begin{lemma}\label{lem:density}
Let $a\ge3$, $t\ge2$, $1\le s\le\min\{t,a-2\}$, and
$M\ge(t+1)(a-1)$. Put $r=a-2$, $q=a+s-1$, and
\[
 D=\binom Mr,\qquad E=\binom Mq,\qquad
 \sigma_s^2=
 \left(\frac{(M-s)_{s-1}}{(M)_{s-1}}\right)^2
 \frac{(r)_s(q)_s}{(M-r)_s(M-q)_s},
\]
If $\mathcal U\subseteq\binom{[M]}r$ and
$\mathcal V\subseteq\binom{[M]}q$ are cross $s$-intersecting, then
their densities $u=|\mathcal U|/D$ and $v=|\mathcal V|/E$ satisfy
\begin{equation}\label{eq:cross-s-density}
 uv\le\sigma_s^2(1-u)(1-v).
\end{equation}
In particular, writing $h=|\mathcal U|$ and $\delta=1-h/D$, we have
\begin{equation}\label{eq:density-missing-factor}
 h|\mathcal V|\le\delta DE\sigma_s^2
 =\delta\,
 \frac{((M-s+1)_s)^2}{(M-r)_s(M-q)_s}
 \binom{M-s}{r-s}\binom{M-s}{q-s}.
\end{equation}
\end{lemma}
\begin{proof}
Equip $\binom{[M]}r$ and $\binom{[M]}q$ with their uniform probability
measures. For $0\le i<s$, let $T_i$ map a function on $\binom{[M]}q$
to its average over the $q$-sets meeting a given $r$-set in exactly
$i$ points.
The hypotheses imply $M\ge r+q$, so these averages are defined.
Both $T_i$ and its adjoint preserve constants.

We construct a signed averaging operator supported on pairs whose
intersection has size less than $s$.
Set
\[
 \mu_\ell=\frac{\binom r\ell\binom q\ell}{\binom M\ell},
 \qquad
 w_i=\sum_{\ell=i}^{s-1}(-1)^{\ell-i}\binom\ell i\mu_\ell,
 \qquad W=\sum_{i=0}^{s-1}w_iT_i.
\]
Binomial inversion gives
\begin{equation}\label{eq:density-moments}
 \sum_{i=0}^{s-1}w_i\binom i\ell=\mu_\ell
 \quad(0\le\ell<s);
\end{equation}
the sum is zero for $\ell\ge s$. In particular $\sum_iw_i=1$,
so $W$ and $W^*$ preserve constants. The weights need not be positive;
we use that $W$ is supported on pairs with intersection smaller
than $s$.

To compute its action in degree $j$, define, for $0\le j\le r$, 
\[
 f_j^{(b)}(A)=\prod_{z=1}^j
      (\mathbf1_{z\in A}-\mathbf1_{j+z\in A}).
\]
By symmetry under interchanging the $j$ pairs $\{z,j+z\}$ and
permuting the remaining points, $T_if_j^{(q)}$ is a multiple of
$f_j^{(r)}$. Evaluate at an $r$-set
containing $[j]$ and avoiding $\{j+1,\ldots,2j\}$ to obtain
$T_if_j^{(q)}=c_j(i)f_j^{(r)}$, where
\begin{equation}\label{eq:density-coefficient}
 c_j(i)=\sum_{z=0}^j(-1)^{j-z}\binom jz
       \frac{(i)_z(q-i)_{j-z}}{(r)_z(M-r)_{j-z}}.
\end{equation}
We use the expansion
\begin{equation}\label{eq:density-newton}
 c_j(i)=(-1)^j\frac{(q)_j}{(M-r)_j}
   \sum_{\ell=0}^j(-1)^\ell\binom j\ell
       \frac{(M+1-j)_\ell(i)_\ell}{(r)_\ell(q)_\ell}.
\end{equation}
To verify this, let $\Delta$ denote the forward difference in $i$.
Using
\[
 \Delta\bigl((i)_z(q-i)_b\bigr)
 =z(i)_{z-1}(q-i-1)_b-b(i)_z(q-i-1)_{b-1},
\]
with terms of zero coefficient omitted, combine adjacent terms in
\eqref{eq:density-coefficient} using
\[
 \frac1{(r)_{z+1}(M-r)_{j-1-z}}
 +\frac1{(r)_z(M-r)_{j-z}}
 =\frac{M+1-j}{r(M-r)(r-1)_z(M-r-1)_{j-1-z}}.
\]
Writing the parameters as superscripts gives
$\Delta c_j^{M,r,q}(i)
=\frac{j(M+1-j)}{r(M-r)}c_{j-1}^{M-2,r-1,q-1}(i)$.
Iterating and evaluating at zero yields
\[
 \Delta^\ell c_j^{M,r,q}(0)
 =(-1)^{j-\ell}\frac{(j)_\ell(M+1-j)_\ell(q)_j}
 {(r)_\ell(q)_\ell(M-r)_j}.
\]
Newton's formula proves \eqref{eq:density-newton}.

Substitute \eqref{eq:density-moments} into
\eqref{eq:density-newton}. The coefficient for $W$ is
\[
 c_j(W)=(-1)^j\frac{(q)_j}{(M-r)_j}
  \sum_{\ell=0}^{\min\{j,s-1\}}(-1)^\ell\binom j\ell
      \frac{(M+1-j)_\ell}{(M)_\ell}.
\]
For $0\le b<j$, induction on $b$ gives the finite identity
\[
 \sum_{\ell=0}^{b}(-1)^\ell\binom j\ell
        \frac{(M+1-j)_\ell}{(M)_\ell}
 =(-1)^b\binom{j-1}b\frac{(M-j)_b}{(M)_b}.
\]
Indeed, subtracting its right side at $b-1$ from that at $b$
recovers the $b$-th summand. The full sum through $b=j$ is zero.
Consequently $c_j(W)=0$ for $1\le j<s$, whereas
\begin{equation}\label{eq:density-operator-coefficient}
 c_j(W)=(-1)^{j+s-1}\binom{j-1}{s-1}
  \frac{(M-j)_{s-1}}{(M)_{s-1}}
  \frac{(q)_j}{(M-r)_j}\qquad(j\ge s).
\end{equation}

To recover the singular values, let $V_j$ be spanned by
$f_J(A)=\mathbf1_{J\subseteq A}$, $|J|\le j$, on the $r$-layer.
For a permutation-invariant self-adjoint operator $K$,
$(Kf_J)(A)$ depends only on $|A\cap J|$. Expand it in
$\binom{|A\cap J|}{\ell}=\sum_{I\subseteq J,\,|I|=\ell}f_I(A)$.
The top coefficient depends only on $|J|$. Thus $K$ preserves $V_j$
and acts as a scalar on $V_j\cap V_{j-1}^{\perp}$.
These spaces, together with the constants, exhaust the function space. The function $f_j^{(r)}$ lies
in this space, since swapping a pair untouched by $I$ proves
orthogonality to $f_I$ for $|I|<j$.

Here $f_j^{(r)}$ and $f_j^{(q)}$ are nonzero since
$j\le r\le\min\{q,M-q\}$. On either layer,
$\|f_j^{(b)}\|^2=2^j\binom{M-2j}{b-j}/\binom Mb$.
By the same symmetry argument as for $W$, $W^*f_j^{(r)}$ is a
multiple of $f_j^{(q)}$. Taking the inner product with $f_j^{(q)}$
and using $Wf_j^{(q)}=c_j(W)f_j^{(r)}$ shows that $WW^*$ has eigenvalue
$c_j(W)^2\|f_j^{(r)}\|^2/\|f_j^{(q)}\|^2$ in degree $j$.
By \eqref{eq:density-operator-coefficient}, all nonzero
nonconstant singular values are therefore among
\begin{equation}\label{eq:density-singular-values}
 \sigma_j=\binom{j-1}{s-1}
   \frac{(M-j)_{s-1}}{(M)_{s-1}}
   \sqrt{\frac{(r)_j(q)_j}{(M-r)_j(M-q)_j}}
       \qquad(s\le j\le r).
\end{equation}
The singular values in degrees $1,\ldots,s-1$ are zero.
This calculation of the singular values uses only
$1\le s\le r\le q$ and $M\ge r+q$.

The largest nonconstant singular value occurs in degree $s$.
Indeed, for $s\le j<r$, consecutive terms in
\eqref{eq:density-singular-values} satisfy
$\frac{\sigma_{j+1}}{\sigma_j} =\frac j{j-s+1}\frac{M-j-s+1}{M-j} \sqrt{\frac{(r-j)(q-j)}{(M-r-j)(M-q-j)}}$.
Each factor is nonincreasing in $j$. For the square-root factor,
pair $r-j$ with $M-q-j$ and $q-j$ with $M-r-j$, using
$M\ge r+q$. It is enough to check the ratio at $j=s$.
For $s=1$ the ratio is at most one. For $s\ge2$, discard
the middle factor and square the remaining expression. Its denominator
is increasing in $M$, so we may substitute
$M=(s+1)(a-1)$, since $t\ge s$. The denominator minus the numerator
is then
$s^2\bigl((s-2)a+4-s\bigr)+s(a-3)>0$.
Thus the norm of $W$ between the mean-zero subspaces is $\sigma_s$.
If $s=r$, there is no consecutive ratio to check.

Finally, put $f=\mathbf1_{\mathcal U}-u$ and
$g=\mathbf1_{\mathcal V}-v$. Since $\mathcal U$ and $\mathcal V$
are cross $s$-intersecting, the support of $W$ gives
$\langle f,Wg\rangle=-uv$. Hence
$uv\le\sigma_s\sqrt{u(1-u)v(1-v)}$.
Squaring and cancelling $uv$ proves \eqref{eq:cross-s-density}
when $uv>0$; otherwise it is immediate. The first inequality in
\eqref{eq:density-missing-factor} follows by discarding $1-v$.
For the binomial form, use
$(M)_s(M-s)_{s-1}/(M)_{s-1}=(M-s+1)_s$ and
$\binom Mr(r)_s=(M)_s\binom{M-s}{r-s}$, and similarly for $q$.
\end{proof}

\subsection{Weighted summation and partial stars}\label{sec:weighted-critical}
\begin{lemma}\label{lem:scalar}
Let $t\ge2$, $A\ge2$, and $M\ge(t+1)A$ be integers. Put
\[
 D=\binom M{A-1},\qquad S=\binom{M+1}A,\qquad
 b=(t+2)\binom{M-1}{A-1}+\binom{M-1}{A-2}.
\]
For $1\le s\le\min(t,A-1)$ define
\[
 \sigma_s^2=
 \left(\frac{(M-s)_{s-1}}{(M)_{s-1}}\right)^2
 \frac{(A-1)_s(A+s)_s}
      {(M-A+1)_s(M-A-s)_s},
 \qquad
 T=\sum_{s=1}^{\min(t,A-1)}
       \binom{t+1}{s+1}D\binom M{A+s}\sigma_s^2.
\]
Then
\[
 (1+D^{-1})T<\frac{t}{t+1}(e-2)b^2<\frac34b^2.
\]
If $T_1$ denotes the first summand of $T$, then
$T_1/b^2\le t/(2(t+1))$. In particular, for every integer $h\ge3$,
\[
 (1+h^{-1})(1+D^{-1})T<b^2.
\]

\end{lemma}
\begin{proof}
Write $C=M-A$, $x=A/C\le1/t$, and
$K_s=\frac{D\binom M{A+s}\sigma_s^2}{S^2}$.
Cancellation gives
$K_1=\frac{A(A-1)C}{(M+1)^2(C-1)}$.
For $1\le s<\min(t,A-1)$, direct cancellation also gives
\[
 \frac{K_{s+1}}{K_s}
 =\frac{A-s-1}{C-2s-1}
   \frac{(C-s)^2}{(C-s+1)(C-2s)}
   \left(\frac{(M-2s+1)(M-2s)}
               {(M-s+1)(M-s)}\right)^2.
\]
All denominators are positive. Since $C\ge2A$,
$\frac{A-s-1}{C-2s-1}\le\frac AC=x$;
indeed, this is equivalent to $C(s+1)\ge A(2s+1)$.
Moreover $s\le t-1$ gives
$M-2s\ge(s+2)A-2s\ge2A$.
Consequently
$\frac{C-s}{C-2s} =1+\frac{s}{C-2s} \le1+\frac{2s}{M-2s} \le\left(\frac{M-s}{M-2s}\right)^2$.
The second and third factors in the displayed ratio for $K_{s+1}/K_s$
therefore have product at most one: use this inequality, together with
$(C-s)/(C-s+1)<1$ and $(M-2s+1)/(M-s+1)<1$.
Thus $K_{s+1}\le xK_s$, and hence
\[
 \frac{T}{S^2}
 \le K_1\sum_{s=1}^t\binom{t+1}{s+1}x^{s-1}
 =K_1\frac{(1+x)^{t+1}-1-(t+1)x}{x^2}.
\]
We have
\[
 \frac bS=\frac A{M+1}L,
 \qquad L=t+2-(t+1)\frac{A-1}{M}.
\]
Thus, with $F_t(x)=(1+x)^{t+1}-1-(t+1)x$,
$\frac{(1+D^{-1})T}{b^2} \le (1+D^{-1})\frac{A-1}{A}\frac C{C-1} \frac{F_t(x)}{x^2L^2}$.
The first three factors have product at most one. Indeed,
$D\ge M\ge3A\ge2(A-1)$, and
$\left(\frac{A-1}{A}\frac C{C-1}\right)^{-1}-1 =\frac{C-A}{(A-1)C}\ge\frac1{2(A-1)}\ge\frac1D$.
Also $L\ge t+1$, since $M\ge(t+1)A$.
The polynomial $F_t(x)/x^2$ has nonnegative coefficients and is
increasing for $x\ge0$. As $x\le1/t$, we obtain
\begin{align*}
 \frac{(1+D^{-1})T}{b^2}
 &\le\frac{t^2}{(t+1)^2}
      \left[\left(1+\frac1t\right)^{t+1}-2-\frac1t\right]\\
 &=\sum_{j=2}^{t+1}
       \frac{t^2}{(t+1)^2}\binom{t+1}j t^{-j}
 <\frac{t}{t+1}\sum_{j=2}^{\infty}\frac1{j!}
 =\frac{t}{t+1}(e-2)<\frac34.
\end{align*}
For the strict termwise comparison, use
$(t+1)_j/t^j\le(t+1)/t$ for $j\ge2$. Keeping only the
$j=2$ term gives $T_1/b^2\le t/(2(t+1))$.
\end{proof}

\begin{proposition}\label{prop:partial}
Assume $a\ge3$, $t\ge2$, $N\ge(t+1)(a-1)+2$. Let
$\mathcal H=(R\vee\mathcal P)\cup\{X\}$, with $|R|=2$, $|X|=a$,
$R\cap X=\varnothing$ and a nonempty
$\mathcal P\subseteq\binom{[N]\setminus R}{a-2}$.
Then $|\mathcal H|W_t(\mathcal H^\perp)<B^2$ if
$|\mathcal P|\ge3$, or $a=3$, or $a=4$ and $|\mathcal P|=2$.
No shifting assumption on $\mathcal P$ is needed.
\end{proposition}
\begin{proof}
Write $A=a-1$, $M=N-2$, $U=[N]\setminus R$, and $h=|\mathcal P|$.
Use $D=\binom M{A-1}$, $T$ and $\sigma_s$ as in
Lemma~\ref{lem:scalar}. We first derive a bound valid for all $h\ge1$. Put
\[
 \mathcal H=(R\vee\mathcal P)\cup\{X\},\qquad
 \mathcal G=\mathcal H^\perp,\qquad
 \mathcal Z=\{Z\in\tbinom Ua:Z\cap P\ne\varnothing
                          \text{ for every }P\in\mathcal P\}.
\]
Put
$Q=\#\{A\in\binom{[N]}a:A\cap R\ne\varnothing, A\cap X\ne\varnothing\} +(t-1)\#\{Y\in\binom U{a-1}:|Y\cap X|\ge2\}$.
By splitting an $(a+d)$-clique according to its $R$-trace, we have
\[
 W_t(\mathcal G)\le
 Q+\sum_{s=1}^{\min(t,a-2)}
       \binom{t+1}{s+1}z_{s-1},
 \qquad z_{s-1}=c_{a+s-1}(\mathcal Z).
\]
The cliques disjoint from $R$ contribute at most $z_d$. Those
meeting $R$ in one point contribute at most $2z_{d-1}$ for $d\ge1$,
and those containing $R$ contribute at most $z_{d-2}$ for $d\ge2$.
Here we omit the additional restrictions imposed by $X$. The cliques
meeting $R$ when $d=0$, and those containing $R$ when $d=1$, give the
first and second terms of $Q$, respectively, including their weights.
The resulting coefficients are
$\binom{t-1}{s-1}+2\binom{t-1}s+\binom{t-1}{s+1}
 =\binom{t+1}{s+1}$.
Every $(a+s-1)$-clique of $\mathcal Z$ meets each member of
$\mathcal P$ in at least $s$ points, since otherwise it would have an
$a$-subset disjoint from that member. In particular the residual
clique family is empty for $s>a-2$.

By Lemma~\ref{lem:density}, for each $s$ in the sum we have
$h z_{s-1}\le(1-h/D) D\binom M{A+s}\sigma_s^2$,
where we have dropped a further factor $1-z_{s-1}/\binom M{A+s}\le1$.
Writing $u=h/D$, we obtain
$(h+1)W_t(\mathcal G) \le(h+1)Q+(1+h^{-1})(1-u)T$.
For the full $R$-star together with $X$, we have $h=D$ and the
companion has weight $Q$. Lemma~\ref{lem:endpoint} gives
$(D+1)Q<B^2$. Hence
$(h+1)W_t(\mathcal G) <\frac{h+1}{D+1}B^2 +\frac{D-h}{D+1}(1+h^{-1})(1+D^{-1})T $.
It suffices to show $(1+h^{-1})(1+D^{-1})T\le B^2$, since
the coefficients $(h+1)/(D+1)$ and $(D-h)/(D+1)$ sum to one.
For $h\ge3$, this follows from Lemma~\ref{lem:scalar}.

\emph{The case $a=4$, $h=2$.}
Retain the factor
\[
 f=\frac{A-1}{A}\frac{C}{C-1},\qquad C=M-A,
\]
from the proof of Lemma~\ref{lem:scalar}. It gives
$\frac{(1+h^{-1})(1+D^{-1})T}{b^2} <(1+h^{-1})(1+D^{-1})f(e-2)$.
Here $A=3$, $C\ge6$ and $D=\binom M2\ge36$, whence
$(1+h^{-1})(1+D^{-1})f \le\frac32\cdot\frac{37}{36}\cdot\frac23\cdot\frac65 =\frac{37}{30}$.
The required ratio is therefore less than
$\frac{37}{30}\cdot\frac34=\frac{37}{40}<1$, so the same
product estimate applies.

\emph{The case $a=3$.}
Here $A=2$, so only $s=1$ occurs. Direct cancellation gives
\[
 D=M,\qquad
 T=\frac{t(t+1)M^2(M-2)}{4(M-3)},\qquad
 b=(t+2)(M-1)+1.
\]
Here $M\ge2(t+1)\ge6$. Put
$L=(t+2)-(t+1)/M\ge t+3/2$, so that $b=ML$. Then
\[
 \frac{2(1+M^{-1})T}{b^2}
 =\left(1+\frac1M\right)\frac{M-2}{M-3}
      \frac{t(t+1)}{2L^2}
 \le\frac76\cdot\frac43
       \frac{t(t+1)}{2(t+3/2)^2}
 <\frac79<1.
\]
Since $1+h^{-1}\le2$ for every $h\ge1$, the same product estimate
applies to every nonempty $\mathcal P$ when $a=3$.
\end{proof}

\subsection{Two larger cores}

\begin{theorem}\label{thm:two-large-cores}
Let $a\ge4$, $t\ge2$, and
\[
 N\ge (t+1)(a-1)+2.
\]
Suppose that $\mathcal H\subseteq\binom{[N]}a$ is the disjoint union
of two nonempty parts $\mathcal H_R$ and $\mathcal H_S$.
Let $R$ and $S$ be disjoint sets of size at least three with
$|R|+|S|\ge a+2$, and suppose that every member of $\mathcal H_R$
contains $R$ and every member of $\mathcal H_S$ contains $S$.
Let $\mathcal G\subseteq\binom{[N]}a$ be cross-intersecting with
$\mathcal H$. Then
\[
 |\mathcal H|W_t(\mathcal G)
 <\left((t+2)\binom{N-3}{a-2}+\binom{N-3}{a-3}\right)^2.
\]
\end{theorem}

\begin{proof}
By replacing $R$ and $S$ by subsets and relabelling the two parts
if necessary, we may assume that their sizes $r,s$ satisfy
$3\le r\le s\le a-1$ and $r+s=a+2$.
Put
$A=a-1,\quad M=N-1,\quad S_0=\binom M A,\quad p=\frac AM\le\frac1{t+1},\quad z=p(1-p),\quad w=\frac p{1-p}$.
Write $b=(t+2)\binom{N-3}{a-2}+\binom{N-3}{a-3}$.
Since
$\frac{\binom{N-3}{a-2}}{S_0} =\frac{(a-1)(N-a)}{(N-1)(N-2)}\ge z$,
we have $b\ge(t+2)zS_0$.

For the core $R$ of size $r$, put $h_R=|\mathcal H_R|$,
$d_r=\binom{N-r}{a-r}$, and $T=t+r-1$.
For $0\le d\le t-1$ partition the $(a+d)$-cliques of $\mathcal G$
by their trace $U$ on $R$, writing $u=|U|$.
If $u\ge d+1$, the trivial bound for that section, multiplied by
$h_R\le d_r$, is
$d_r\binom{N-r}{a+d-u}$.
If $u\le d$, put $j=d+1-u\ge1$. The family
$\{H\setminus R:H\in\mathcal H_R\}$ and the family obtained by
deleting $U$ from the cliques in this section are cross $j$-intersecting
on $[N]\setminus R$. Their uniformities are $a-r$ and
$a+d-u=a+j-1$, respectively.
If $j>a-r$, the section is empty. Otherwise, when
$N-r\ge(j+1)a$, inequality~\eqref{eq:unequal-rank-ekr} gives
\begin{equation}
 h_R\,|\text{section}|
 \le \binom{N-r-j}{a-r-j}\binom{N-r-j}{a-1}.             \label{eq:formula-l1}
\end{equation}
All these applications are valid except possibly $j=t$.
Indeed, for $j\le\min\{t-1,a-r\}$,
$N-r-(j+1)a\ge(t-j)a-t-r+1\ge0$.
For the last inequality, if $a-r\ge t-1$, use $j\le t-1$ and
$a\ge t+r-1$. Otherwise put $v=t+r-a\ge2$ and use $j\le a-r$;
the resulting lower bound is $(v-1)(a-1)\ge0$.
The exceptional case $j=t$ can only be $d=t-1,u=0$.
Embed its two families in an ambient set of size $n=N+t-1$.
Now $n\ge(t+1)a$, so we can apply \eqref{eq:unequal-rank-ekr}
with larger uniformity $a+t-1$ and obtain
\begin{equation}
 h_R\,|\text{exceptional section}|
 \le \binom{N-1}{a-r-t}S_0\le w^{r+t-1}S_0^2.          \label{eq:formula-l2}
\end{equation}
The section is empty if $a-r<t$.

Put $i=r+d-u$. Summing over $d$ and $u$ gives
$\sum_{d-u=i-r}\binom{t-1}d\binom ru=\binom Ti$.
Summing the preceding bounds therefore gives
\begin{align}
 h_R W_t(\mathcal G)
 &\le d_r\sum_{i=0}^{r-1}\binom Ti\binom{N-r}{a-r+i}\\
 &\quad+\sum_{i=r}^{a-1}\binom Ti
      \binom{N-i-1}{a-i-1}\binom{N-i-1}{a-1}
      +w^{r+t-1}S_0^2.                                  \label{eq:formula-l3}
\end{align}
The second sum may include the formal term for $j=t$; adding this
nonnegative term only weakens the bound.

Dividing the binomial products in \eqref{eq:formula-l3} by
$S_0^2$, we have
\begin{equation}
 \frac{d_r\binom{N-r}{a-r+i}}{S_0^2}
 \le p^{2(r-1-i)}z^i \quad(0\le i<r),\qquad
 \frac{\binom{N-i-1}{a-i-1}\binom{N-i-1}{a-1}}{S_0^2}
 \le z^i \quad(r\le i\le a-1).                            \label{eq:formula-l4}
\end{equation}
For the first, the exact ratio is
$\frac{(A)_{r-1}(A)_{r-i-1}(M-A)_i}{(M)_{r-1}^2}$.
For each $0\le v\le r-2$, $A/M=p\le1/2$ gives
$1-v/A\le(1-v/M)^2$; the factors of $(A)_{r-1}/A^{r-1}$ are therefore bounded by those of
$((M)_{r-1}/M^{r-1})^2$, while the remaining normalized falling
factorials are at most one. This proves the first estimate.
The ratio in the second estimate is
$(A)_i(M-A)_i/(M)_i^2$, which is at most $p^i(1-p)^i=z^i$.

Define
\[
 G_r=\sum_{i=0}^{r-1}\binom{t+r-1}i p^{2(r-1-i)}z^i
       +\sum_{i=r}^{t+r-1}\binom{t+r-1}i z^i.
\]
Then \eqref{eq:formula-l3}--\eqref{eq:formula-l4} yield
\begin{equation}
 h_R W_t(\mathcal G)\le (G_r+w^{r+t-1})S_0^2.           \label{eq:formula-l5}
\end{equation}

Suppose $a\ge5$ and write $G_r=A_r+B_r$ for its two displayed sums.
Pascal's identity gives $G_{r+1}=pA_r+(1+z)B_r$: a term of $A_r$
has factors $p^2$ and $z$ in the next sum, whereas a term of $B_r$
has factors $1$ and $z$.
The first term of $B_r$ divided by the final term of $A_r$ is $tz/r$;
ratios of successive terms within $B_r$ are at most
$(t-1)z/(r+1)<1/(r+1)$. Since $tz<1$, a geometric series gives
$B_r\le \frac{r+1}{r^2}A_r\le\frac{A_r}{r-1}$.
Hence $G_{r+1}-G_r
=-(1-p)A_r+zB_r<0$ for $r\ge2$.
Here $r\ge3$ and $s\ge4$, so \eqref{eq:formula-l5} for the two parts gives
\begin{equation}
 \frac{|\mathcal H|W_t(\mathcal G)}{S_0^2}
 \le G_3+G_4+w^{t+2}+w^{t+3}.                          \label{eq:formula-l6}
\end{equation}

To bound the right side by $(t+2)^2z^2$, first take $t=2$
and divide by $16z^2$. Every resulting summand is nondecreasing
in $p\in(0,1/3]$: the terms with powers $z^i$, $i\ge2$, give
nonnegative powers of $z$ after division. Each remaining term is a
positive constant times a positive power of $p$ divided by a positive
power of $1-p$. Evaluation at $p=1/3$ gives exactly
\begin{equation}
 \frac{G_3+G_4+w^4+w^5}{16z^2}
 \le\frac{1398139}{1492992}<1.                          \label{eq:formula-l7}
\end{equation}

Suppose $t\ge3$. Use $p/(1-p)\le1/t$,
$p\le1/(t+1)$, $z\le t/(t+1)^2$, and $tz<1$.
Expanding the first few terms and bounding the remaining sums by
the exponential series gives
\begin{align*}
 \frac{G_3}{z^2}
 &\le\frac1{t^2}+\frac{t+2}{t}+(t+2)(t+1)(e-2),\\
 \frac{G_4}{z^2}
 &\le\frac1{t^2(t+1)^2}+\frac{t+3}{t(t+1)^2}
       +\frac{(t+3)(t+2)}{2(t+1)^2}
       +\frac{(t+3)(t+2)t}{t+1}(e-5/2).
\end{align*}
For example, after its first two factors,
$\binom{t+2}i\le(t+2)(t+1)t^{i-2}/i!$; the corresponding
three-factor estimate gives the second line's tail.
Substituting $e<11/4$ and simplifying yields
\begin{equation}
 \frac{G_3+G_4}{(t+2)^2z^2}
 <1-\frac{(t-2)(3t+2)}{4t^2(t+2)}.                    \label{eq:formula-l8}
\end{equation}
The remaining correction satisfies
\begin{equation}
 \frac{w^{t+2}+w^{t+3}}{(t+2)^2z^2}
 \le\frac{(t+1)^5}{(t+2)^2t^{t+5}}
 <\frac5{(t+2)^2t^3}<\frac1{4(t+2)}
 \le\frac{(t-2)(3t+2)}{4t^2(t+2)}.                    \label{eq:formula-l9}
\end{equation}
This proves the desired strict inequality in \eqref{eq:formula-l6}.

It remains to consider $a=4$.
Here $r=s=3$. Put $m=N-3$ and $T=t+2$; then
$m\ge3t+2=3T-4\ge2T$ and $T\ge4$.
The tail rank is one, so only $j=1$ occurs in \eqref{eq:formula-l1}, and there is
no exceptional term. Formula \eqref{eq:formula-l3} becomes, for each part,
$h_R W_t(\mathcal G)\le F:= m^2+Tm\binom m2+\binom T2m\binom m3 +\binom T3\binom{m-1}3$.
Since $b=T\binom m2+m$, direct termwise estimates give
\begin{align*}
 \frac{2F}{T^2\binom m2^2}
 &=\frac8{T^2(m-1)^2}+\frac4{T(m-1)}
   +\frac23\frac{T-1}{T}\frac{m-2}{m-1}\\
 &\quad+\frac29\frac{(T-1)(T-2)}T
                  \frac{(m-2)(m-3)}{m^2(m-1)}\\
 &<\frac1{98}+\frac17+\frac23+\frac19
   =\frac{821}{882}<1.
\end{align*}
Thus $|\mathcal H|W_t(\mathcal G)\le2F<b^2$ in this case as well.
\end{proof}

\medskip\label{sec:completion}
\begin{proof}[Proof of Theorem~\ref{thm:critical}]
Enlarge $\calG$ to $\calH^\perp$; this can only increase $W_t(\calG)$.
Partition $\calH$ according to \eqref{eq:critical-core-structure}.
If $a\ge4$ and both cores have size at least three, Theorem~\ref{thm:two-large-cores} gives $|\calH|W_t(\calG)<b^2\le B^2$.
For $a=3$, two cores of size three are two disjoint triples; shrinking
one of these cores to size two places the pair in \eqref{eq:partial-star-form}.
If one of the two cores has size two, the partition already has
the form \eqref{eq:partial-star-form}.

Write $h_P=|\mathcal P|$ and $K=\bigcap\mathcal P$.
Proposition~\ref{prop:partial} covers $h_P\ge3$, every $h_P\ge1$
when $a=3$, and $h_P=2$ when $a=4$.
It remains to treat $a\ge4$ with $K\ne\varnothing$, and then
$a\ge5$, $h_P=2$, $K=\varnothing$.

If $K\ne\varnothing$, nontriviality implies $K\cap X=\varnothing$.
Choose one point $x\in K$. The two parts in \eqref{eq:partial-star-form} have disjoint cores
$R\cup\{x\}$ and any $(a-1)$-subset of $X$. Their sizes are
$3,a-1$, so Theorem~\ref{thm:two-large-cores} applies. This includes $h_P=1$.

Finally, suppose $a\ge5$, $h_P=2$, and $K=\varnothing$.
Write $\mathcal P=\{P_1,P_2\}$; then $P_1,P_2$ are disjoint
$(a-2)$-sets. We claim that one is contained in $X$.
Otherwise choose $u\in P_1\setminus X$ and $v\in P_2\setminus X$,
and interchange the subscripts so that $u<v$.
The left shift $s_{uv}$ fixes $X$ and $R\cup P_1$, but changes
$R\cup P_2$ to $R\cup(P_2\setminus\{v\})\cup\{u\}$.
The new member is absent from the original family, since $a-2\ge2$
and the tails are disjoint. The intersection of the two star members
after shifting is $R\cup\{u\}$, which is disjoint from $X$.
Thus the shift changes $\calH$ and leaves it nontrivial, contrary
to the assumption that $\calH$ is critically shifted.
We may therefore assume $P_1\subseteq X$.
Repartition the family as
$\{R\cup P_1,X\}\ \dot\cup\ \{R\cup P_2\}$.
The first part has core $P_1$, of size $a-2\ge3$, disjoint from
the singleton in the second part. Choose disjoint cores of sizes
$3,a-1$ within these two sets. Theorem~\ref{thm:two-large-cores} applies
once more, completing the weighted estimate.

Now let $(\mathcal A,\mathcal B)$ be a saturated terminal pair of
nontrivial cross $t$-intersecting families, and suppose that a
simultaneous shift makes $\mathcal A$ trivial. The argument in
Case~2 of the proof of Theorem~\ref{thm:simplified-weighted-reduction}
uses only saturation and this shift, and gives
\[
 C=\bigcap\mathcal A,\qquad |C|=t-1,\qquad
 \mathcal A=C\vee\mathcal H=L_C(\mathcal H),
\]
where $\mathcal H$ is nontrivial and admits a trivializing shift.
By saturation and Lemma~\ref{lem:clique-lift-duality},
\[
 \mathcal B=L_C(\mathcal G),\qquad
 \mathcal G=\mathcal H^\perp,\qquad
 \mathcal H=\mathcal G^\perp,
\]
and $\mathcal G$ is nontrivial. Thus
$|\mathcal A||\mathcal B|=|\mathcal H|W_t(\mathcal G)$.
To prove that $\mathcal H$ is critically shifted, first note that
$(\mathcal H,\mathcal G)$ is terminal for shifts on
$[n]\setminus C$. Such a shift leaves the trace on $C$ unchanged
and restricts on the base layers to the shifts of $\mathcal H$
and $\mathcal G$. If it changed the residual pair and left both
intersections empty, it would change the original pair and leave
both families nontrivial: the shifted base layer in each family
would have intersection exactly $C$.
By \cite[Lemma 3.2]{WZ}, in a saturated terminal pair of ordinary
cross-intersecting $a$-uniform families on $N\ge2a$ points, at most
one family admits a shift that makes it trivial. Here
$N\ge(t+1)(a-1)+2\ge2a$. Since $\calH$ admits such a shift,
no shift can make $\calG$ trivial. It follows that every shift
that changes $\calH$ makes it trivial. Thus $\calH$ is critically
shifted, and the weighted estimate proves the assertion.
\end{proof}

\section{Shifted families and kernel estimates}\label{sec:shifted}
Throughout this section we assume
\[
 t\ge2,\qquad k\ge t+2,\qquad n\ge(t+1)(k-t+1),
\]
and write $\Phi=m_t(n,k)^2$. Note that
$n-2k\ge(t-1)(k-t-1)\ge1$. By Theorem~\ref{thm:critical}, it
suffices to consider shifted extremal pairs. We use generator
exchanges to reduce to ordinary kernels with a universal vertex.
Proposition~\ref{prop:universal-vertex-closure} bounds their products
using the estimates proved later in this section.

\subsection{Terminal configurations}
\begin{proposition}\label{prop:shifted-bound}
Under the standing assumptions, let $\mathcal A$ and $\mathcal B$
be nonempty, nontrivial, shifted cross $t$-intersecting families.
Then $|\mathcal A||\mathcal B|\le\Phi$.
\end{proposition}

The following terminal configuration will be used in the proof.
Let $(\mathcal A,\mathcal B)$ be a saturated pair of nonempty,
nontrivial, shifted cross $t$-intersecting families with common
maximal canonical support $[s]$, where $s\ge t+4$. Suppose that every canonical generator of $\mathcal A$ containing $s$ has size $t+1$, every canonical
generator of $\mathcal B$ containing $s$ has size $s-1$, and all generators of $\mathcal B$
have size at least $t+2$. Let $\mathcal B_0$ be generated by the
generators of $\mathcal B$ avoiding $s$, and suppose that it is trivial.
By Lemma~\ref{lem:terminal-deletion}, $\mathcal B_0$ is nonempty;
it is shifted, so
\[
 C=\bigcap\mathcal B_0=[c],\qquad c\ge t.
\]
The active layer of $\mathcal B$ also gives $s\le k+1$. We use these assumptions
and notation in the next two lemmas.

\begin{lemma}\label{lem:exact-terminal-layers}
Under the preceding assumptions,
\[
 \{F\in\mathscr G(\mathcal B):s\in F\}
       =\{[s]\setminus\{p\}:p\in C\},
\]
\begin{equation}
 \{E\in\mathscr G(\mathcal A):s\in E\}
       =\{U\cup\{s\}:U\in\tbinom Ct\}.                     \label{eq:formula-s4}
\end{equation}
\end{lemma}
\begin{proof}
Every $(s-1)$-subset of $[s]$ meets every generator of $\mathcal A$ in at least
$t$ points, so saturation includes its complete star in $\mathcal B$.
For $p\in C$, no retained generator is contained in $[s]\setminus\{p\}$,
and all generators containing $s$ have size $s-1$, so this co-singleton
is canonical. For $p\notin C$, some retained generator avoids $p$ and
is contained in that co-singleton. Hence the co-singleton either omits
$s$ or is not canonical.

An active generator of $\mathcal A$ has the form $U\cup\{s\}$ with $|U|=t$.
It must meet every retained generator of $\mathcal B$ in at least $t$ points,
so $U\subseteq C$. Conversely, every such $U\cup\{s\}$ meets
each retained generator and each co-singleton in at least $t$ points.
Saturation adds its full star, and the absence of generators of
rank at most $t$ makes it canonical.
\end{proof}

For the ordinary sections of this pair, put
\[
 a=k-t+1,\quad N=n-t+1,\quad x=t,\quad
 U=[t+1,n],\quad X=[t+1,s],\quad m=|X|=s-t\le a.
\]
Take the full $[t-1]$-sections of $\mathcal A$ and $\mathcal B$,
and delete $[t-1]$ from them; denote the resulting ordinary
$a$-uniform families by $\mathcal H$ and $\mathcal G$.
The shifted-pair argument in the proof of
Theorem~\ref{thm:simplified-weighted-reduction} applies to any
saturated shifted pair. It shows that $\mathcal H,\mathcal G$ are
shifted, saturated, nontrivial cross-intersecting families and that
their clique lifts recover $\mathcal A,\mathcal B$.
In particular,
\begin{equation}
 |\mathcal A|=W_t(\mathcal H),\qquad
 |\mathcal B|=W_t(\mathcal G),\qquad
 W_t(\mathcal L)=\sum_{d=0}^{t-1}\binom{t-1}d c_{a+d}(\mathcal L).
                                                                    \label{eq:formula-s5}
\end{equation}
For a family $\mathcal L$ and a point $v$, write
\[
 \mathcal L(v)=\{L\setminus\{v\}:L\in\mathcal L,\ v\in L\},
 \qquad
 \mathcal L(\overline v)=\{L\in\mathcal L:v\notin L\}.
\]
For a family $\mathcal Q$ on $U$, write $\mathcal Q^{\perp_r}$ for
its rank-$r$ ordinary cross-intersecting companion. Set
\[
 \mathcal E=\mathcal H(\overline x),\qquad
 \mathcal F=\{F\in\tbinom Ua:X\subseteq F\},\qquad
 \mathcal P=\mathcal E^{\perp_{a-1}}.
\]
All retained generators of $\mathcal B$ contain $[t]$.
By Lemma~\ref{lem:exact-terminal-layers}, the only active generator
that can occur in a member containing $[t-1]$ and omitting $x=t$
is $[s]\setminus\{t\}$. Thus the section omitting $x$ is $\mathcal F$.
Saturation then gives
\begin{equation}
 \mathcal G(\overline x)=\mathcal F,\qquad
 \mathcal H(x)=\mathcal F^{\perp_{a-1}}
             =\{Q\in\tbinom U{a-1}:Q\cap X\ne\varnothing\},
 \qquad \mathcal G(x)=\mathcal P.                              \label{eq:formula-s6}
\end{equation}
Both $\mathcal E$ and $\mathcal P$ are generated by subsets of $X$:
take the traces on $X$ of canonical generators compatible with the
respective sections. Here a kernel on $X$ is an upset in $2^X$, specified by its minimal
members; its rank-$r$ extension on $U$ consists of the $r$-sets
containing one of these members. The blocker of a kernel is the
upset whose minimal members are the minimal transversals of its
minimal generators.
The minimal generators of $\mathcal E$ and $\mathcal P$ have size
at least two. A generator of size at most one in $\mathcal E$ would give a
complete star in $\mathcal A$ through a set of size at most $t$,
which is impossible.
A generator of size at most one in $\mathcal P$ would give a
generator of $\mathcal B$ of rank at most $t+1$, contrary to our assumption.
In the first assertion, any extension of the proposed generator
through $x$ is also in $\mathcal A$: replace $x$ by an unused point of
$U$ and shift back. Their kernels are mutual blockers. Indeed, every minimal transversal
of a kernel on $X$ with no singleton generator has size at most
$m-1\le a-1$: each set $X\setminus\{u\}$ is already a transversal.
Moreover, $|U|\ge2a$, so a prescribed generator disjoint from
a tested set of size at most $a$ can be extended, still disjoint,
to the required uniform rank. Thus the fixed-rank companion
condition is exactly the blocker condition. The rank restrictions
lose no minimal transversal, and saturation gives the mutual blocker relation.

\begin{lemma}
\label{lem:residual-terminal-pair}
Let $R$ be a minimal residual generator of $\mathcal E$ containing
the last point $s$. Then $[t-1]\cup R$ is a canonical generator of $\mathcal A$.
Consequently $|R|=2$, and its possible other points are exactly
$C\setminus[t]$.
\end{lemma}
\begin{proof}
Choose a canonical generator of $\mathcal A$ $E\subseteq[t-1]\cup R$.
Residual minimality forces $R\subseteq E$. If $p\in[t-1]\setminus E$,
shift $s$ to $p$. This produces a complete-star generator contained
in $[t-1]\cup(R\setminus\{s\})$, contradicting residual minimality.
Thus $E=[t-1]\cup R$. The terminal layer of $\mathcal A$ has rank $t+1$, so
$|R|=2$. By \eqref{eq:formula-s4}, the other point lies in
$C\setminus[t]$. Conversely, for every $u\in C\setminus[t]$,
$[t-1]\cup\{u,s\}$ is one of the generators in
\eqref{eq:formula-s4}, and $\{u,s\}$ is minimal in the residual
kernel because it has no generator of size less than two.
\end{proof}

\subsection{The shifted bound and terminal kernels}
\begin{proof}[Proof of Proposition~\ref{prop:shifted-bound}]
Suppose that the maximum product over all nonempty nontrivial cross
$t$-intersecting pairs exceeds $\Phi$. Applying simultaneous shifts
that preserve nontriviality to a maximizing pair
gives a saturated terminal pair $(\mathcal A,\mathcal B)$ with the
same product. Theorem~\ref{thm:critical} forces this pair to be shifted.

Use the canonical generators from Lemma~\ref{lem:canonical-prefix},
and let $s$ be their common maximal support. By
Lemma~\ref{lem:canonical-pairing}, all generator ranks are at least
$t+1$, $s\ge t+2$, and active ranks are paired by
\begin{equation}
 i+j=s+t.                                                     \label{eq:formula-s1}
\end{equation}
If each family has a $(t+1)$-generator, shiftedness supplies the
complete $[t+1]$-star in each family. Every member of both families
then meets $[t+1]$ in at least $t$ points, and Theorem~\ref{thm:transversal}
gives the required bound. Interchanging the two families if necessary,
we may thus assume that
\begin{equation}
 \min_{E\in\mathscr G(\mathcal B)}|E|\ge t+2.                 \label{eq:formula-s2}
\end{equation}

The case $s=t+2$ is excluded by \eqref{eq:formula-s1} and
\eqref{eq:formula-s2}, since every generator of $\mathcal A$ has size at least $t+1$.
When $s=t+3$, condition \eqref{eq:formula-s2} and saturation force
\[
 \mathcal A=\{A:|A\cap[t+3]|\ge t+1\},\qquad
 \mathcal B=\{B:|B\cap[t+3]|\ge t+2\}.
\]
Indeed every $(t+1)$-subset of $[t+3]$ meets every generator of $\mathcal B$ in
at least $t$ points, and every $(t+2)$-subset meets every generator of $\mathcal A$
in at least $t$ points. Lemma~\ref{lem:asymmetric-threshold} bounds
this product by $b_t(n,k)^2$. Hence $s\ge t+4$.

There can be no active pair of ranks
\begin{equation}
 t+2\le i,j\le s-2.                                          \label{eq:formula-s3}
\end{equation}
By Lemma~\ref{lem:terminal-deletion}, both exchanges leave each
family nontrivial. Lemma~\ref{lem:paired-terminal-exchange} gives
the changes in the sizes of the families and the shadow estimates. For $t\ge3$, these contradict
Lemma~\ref{lem:general-individual-exchange}. For $t=2$,
Lemma~\ref{lem:t2-individual-exchange} handles the nonexceptional
ranks. Lemma~\ref{lem:t2-global-middle} and the argument following it
rule out the exceptional ranks, either by improving the product or
by showing it is below $b_2(n,k)^2$.

By \eqref{eq:formula-s1}--\eqref{eq:formula-s3}, every generator of $\mathcal A$
containing $s$ has size $t+1$ and every generator of $\mathcal B$ containing $s$
has size $s-1$. By Lemma~\ref{lem:paired-terminal-exchange},
the retained family $\mathcal B_0$ is cross $t$-intersecting with
the full star through $E\setminus\{s\}$ for each active generator of $\mathcal A$
$E$. Since $|E\setminus\{s\}|=t$ and $n>2k-t$, every member of
$\mathcal B_0$ contains $E\setminus\{s\}$. Thus $\mathcal B_0$ is trivial, giving the terminal configuration
above. Retain its notation
$C=[c]$, $a,N,x,U,X,m$, and $\mathcal E,\mathcal F,\mathcal P$,
with the exact sections in \eqref{eq:formula-s6}.

If $c>t$, the residual last point therefore has the nonempty initial
neighbor set $D=C\setminus[t]$. By shiftedness, every two-element subset of $X$ meeting $D$
belongs to the kernel. Every other minimal generator avoids
$D\cup\{s\}$ and has size at least two. This is precisely the
universal-vertex decomposition treated by
Proposition~\ref{prop:universal-vertex-closure}, including the case
of an empty residual kernel. It remains to consider $c=t$.

Assume $C=T=[t]$. Lemma~\ref{lem:residual-terminal-pair} shows that
$s$ is unused by $\mathcal E$. Let $Y$ be its actual initial support,
of size $y<m$. Since its blocker $\mathcal P$ has no singleton
generator, $\mathcal E$ has no common point and $y\ge3$.
Moreover, $\mathcal B_0$ consists of exactly the members of
$\mathcal B$ containing $T$. Indeed, such a member that contains an
active generator of $\mathcal B$ contains all of $[s]$ and hence a retained
generator as well. Thus $\mathcal B_0=\{T\cup Q:Q\in\mathcal P\}$.
Since its core is $T$, the members of $\mathcal P$ have no common point.

Fix $\mathcal E$ and $\mathcal P$, and allow the outer initial core
$X$ to have any size $z$ with $y\le z\le a$. Define the ordinary
pair by \eqref{eq:formula-s6} with this new $X$, leaving $\mathcal E$ fixed.
Each member of $\mathcal E$ meets $X$ in at least two points, and
each member of $\mathcal F$ contains $X$. Thus the two sections
avoiding $x$ are cross-intersecting; the other section pairs are
cross-intersecting by construction. The pairs remain shifted:
inserting $x$ into a member of $\mathcal E$ leaves a point of $X$,
and inserting it into a member of $\mathcal F$ leaves at least one
point of each generator of $\mathcal E$. Shifts within $U$ are
allowed because $X$ is initial and the sections there are shifted.
The family $\mathcal H$ contains the distinct stars through
$\{x,u\}$, $u\in X$, and a member avoiding $x$. Similarly, the
section $\mathcal P$ has empty intersection and $\mathcal F$ is
nonempty, so $\mathcal G$ is nontrivial as well.
Their lifts are shifted: shifts into $[t-1]$ delete a point of the
outside clique, and the other shifts preserve the clique condition.
Lemma~\ref{lem:clique-lift-duality} gives nontrivial cross
$t$-intersection.

For either lift, membership depends only on the trace on $[t]\cup X$:
permuting points outside this block preserves all the ordinary
sections and hence the clique condition. All canonical generators
are supported on this block. To see this, let $E$ be a complete-star
generator and put $E_0=E\cap([t]\cup X)$. Since $n\ge 2k+1$ and
$t+z\le k+1$, we can extend $E$ to a member having trace exactly
$E_0$. Permuting outside points and then shifting additional points
into the block gives every $k$-set containing $E_0$. Thus $E_0$ is
itself a complete-star generator.

Write $M=N-1$,
\[
 S=\binom M{a-1},\qquad Q=|\mathcal P|,\qquad
 J=\sum_{j=0}^{t-1}\binom t{j+1}c_{a+j}(\mathcal E).
\]
For $d\ge1$, an $(a+d)$-clique of $\mathcal H$ through $x$
is $\{x\}\cup K$ for an $(a+d-1)$-clique $K$ of $\mathcal E$.
Indeed, every $(a-1)$-subset of $K$ meets $X$: otherwise it would
extend in $K$ to an $a$-subset meeting $X$ in at most one point.
There is no clique of $\mathcal F$ of size greater than $a$, as every
member of $\mathcal F$ contains $X$. Therefore the only higher cliques
of $\mathcal G$ are $\{x\}\cup F$ for $F\in\mathcal F$; every
$(a-1)$-subset of $F$ meets each generator of $\mathcal E$, since
$X\subseteq F$ and each such generator has at least two points in $X$.
Summing the clique layers gives
\begin{equation}
 V_z=\left(S-\binom{M-z}{a-1}+J\right)
       \left(Q+t\binom{M-z}{a-z}\right).                     \label{eq:formula-s7}
\end{equation}
Proposition~\ref{prop:unused-core-endpoints} gives
\begin{equation}
 V_m\le\max\{V_y,V_a\}.                                    \label{eq:formula-s8}
\end{equation}
By maximality, the pair at an endpoint attaining the right-hand side
also has maximum product. In particular, it is already saturated.
At $z=y$, the common maximal canonical support is $[t+y]$.
It is contained in this interval by the support argument above.
For the reverse inclusion, take a minimal generator $R$ of
$\mathcal E$ containing $t+y$ and extend it to an $a$-set with
trace $R$ on $Y$. Replacing $t+y$ by an unused point outside $Y$
changes its membership in $\mathcal E$, because $R\setminus\{t+y\}$
does not contain a generator. Appending $[t-1]$ gives two $k$-sets
with the same trace on $[t+y-1]$, only one of which belongs to
$\mathcal A$. The new family $\mathcal B$ has
no $(t+1)$-generator. Indeed, shifting such a generator to contain
$[t-1]$ would give a residual two-element generator of $\mathcal G$.
This is impossible: $\mathcal P$ has no singleton generator and
$\mathcal F$ requires the whole set $X$ of size $y\ge3$.
If $y=3$, the asymmetric threshold case applies. If $y\ge4$,
apply the terminal exchanges again to the endpoint pair. They apply
because this pair is shifted and attains the global maximum.
The terminal layer of $\mathcal A$ has rank $t+1$, so the proof of
Lemma~\ref{lem:residual-terminal-pair} shows that the last point of
$Y$ occurs only in pairs in the residual kernel. Since it occurs in
at least one such pair, shiftedness again gives the universal-vertex
decomposition and Proposition~\ref{prop:universal-vertex-closure}
applies.

It remains to consider $z=a$, where the outer support is $k+1$
and $Y$ is a proper initial subset of the outer residual core.
The members of $\mathcal B$ outside
$\mathcal B_0=\{T\cup Q:Q\in\mathcal P\}$ are exactly the $t$
sets $[k+1]\setminus\{p\}$, $p\in T$.
Apply a second exchange on the actual residual support, with
$q=t+y<k+1$. The family $\mathcal B_0$
has core $T$ and canonical generators supported on $[q]$, since its
membership is determined by the trace there. By
Lemma~\ref{lem:retained-canonical}, these generators are precisely
those retained after the outer exceptions are removed.
The mutual blocker description shows that its actual residual
support equals $Y$. Indeed, for a point $v$ in a minimal generator
$R$ of a clutter on $Y$, the set $(Y\setminus R)\cup\{v\}$ meets
every clutter member. A minimal transversal inside it must retain
$v$, since it must meet $R$. Thus a clutter and its blocker have
exactly the same support.
Every $T\cup\{u\}$ with $u\in X$ is a canonical generator
of $\mathcal A$: its complete star is in the section through $[t]$,
and no canonical generator has size at most $t$.
Conversely, let $E$ be a canonical generator of $\mathcal A$ with
$u=\max E>q$. If $|E|\ge t+2$, then $E-u$ still crosses every
retained generator of $\mathcal B$, since they avoid $u$, and
every outer co-singleton, since $|E-u|\ge t+1$ and
$E-u\subseteq[k+1]$. Saturation would contradict minimality of $E$.
Thus $|E|=t+1$. Now $E-u$ must be contained in every retained
generator of $\mathcal B$, and therefore equals $T$.

Suppose that a generator $E$ of $\mathcal A$ with maximum $q$ has rank
$i\ge t+2$. Shortening it by $q$ still meets every outer B-exception
in at least $t$ points. Minimality and saturation therefore supply
a retained generator $F$ of $\mathcal B$ with $q\in F$ and $|E\cap F|=t$.
As in Lemma~\ref{lem:canonical-pairing}, shifting $q$ into a missing
point proves
\begin{equation}
 E\cup F=[q],\qquad j=|F|=q+t-i,
 \qquad t+2\le i,j\le q-2.                                  \label{eq:formula-s9}
\end{equation}
In particular $q\ge t+4$; when $q=t+3$, such activity is already
impossible.

Let $I,J$ be the two paired inner layers, of sizes $\alpha,\beta$.
With cylinders defined by the exact prefix on $[q]$, the four changes
in the forward and reverse exchanges are
\begin{equation}
 \alpha\binom{n-q}{k-i+1},\quad
 \alpha\binom{n-q}{k-i},\quad
 \beta\binom{n-q}{k-j+1},\quad
 \beta\binom{n-q}{k-j}.                                     \label{eq:formula-s10}
\end{equation}
These are exact counts. Each selected generator of $\mathcal A$
omits a point of $T$, as it cannot contain a retained generator
$T\cup\{u\}$. Its cylinders thus avoid the fixed outer stars.
Each affected cylinder of $\mathcal B$ contains $T$ and therefore
avoids the outer exceptions. An additional prefix point supplies a
retained generator by shifting $q$ to that point. Conversely, an
exact prefix cannot contain another generator, by minimality.
For cross $t$-intersection, any inner opposite generator that
violates the bound with a shortened generator must contain $q$
and meet the original generator in exactly $t$ points. The union
argument of Lemma~\ref{lem:paired-terminal-exchange} then gives
union $[q]$, so that this opposite generator lies in the layer
being deleted. The outer generators cause no obstruction:
shortened generators of $\mathcal A$ have at least $t+1$ points
in $[k+1]$, and shortened retained generators of $\mathcal B$
still contain $T$.

Both deletions leave the corresponding families nontrivial. On the A-side, all
low generators $T\cup\{u\}$, $t<u\le k+1$, survive, so any common
$t$-set after deletion must be $T\subseteq[q-1]$. On the B-side,
the outer exceptions survive, so a common $t$-set must avoid $T$;
shifting $q$ out of a deleted generator also supplies a surviving
generator inside $[q-1]$, so that common set lies in $[q-1]$.
Let $C_0\subseteq[q-1]$ be such a common $t$-set on either side.
Since the original family is nontrivial, a deleted generator $R$
omits some $c_0\in C_0$. For every $v\in[q-1]\setminus R$,
shifting $q$ to $v$ gives a surviving complete-star generator inside
$R-q+v$. Hence $C_0\subseteq R-q+v$, forcing $v=c_0$.
Thus $R=[q]\setminus\{c_0\}$, contrary to
\eqref{eq:formula-s9}.

Finally, let $I'=\{E-q:E\in I\}$ and let $\partial^+I'$ be its
rank-$i$ upper shadow in $[q-1]$. Shiftedness and normalized matching
give
\begin{equation}
 |\partial^+I'|\ge\alpha\frac{q-i}{i},\qquad
 |\mathcal A|\ge\alpha\binom{n-q}{k-i}
       +|\partial^+I'|\binom{n-q+1}{k-i}.                    \label{eq:formula-s11}
\end{equation}
The two counted parts have different prefix sizes and are disjoint.
The same bound holds on the B-side. We may therefore use the
smaller of $i,j$ as the first rank in the algebraic estimate.
Lemma~\ref{lem:general-individual-exchange}
applies with $s=q$ when $t\ge3$. For $t=2$, the numerical
argument is supplied directly by
Lemma~\ref{lem:t2-individual-exchange}; the exceptional balanced
inner ranks are covered by Lemma~\ref{lem:t2-inner-extra-cylinder}.
Its additional cylinders are supplied by the retained generators
$[2]\cup\{u\}$ with $q<u\le k+1$.
One of the two exchanges strictly increases the product, a contradiction.

Consequently every generator of $\mathcal A$ with maximum $q$ has size $t+1$. The proof of
Lemma~\ref{lem:residual-terminal-pair}, now at $q$, shows that every
minimal residual generator using that point has size two. Since $Y$
is the actual support, there is at least one. Its shifted neighbor
set is a nonempty initial interval $D$; all pairs meeting $D$ occur,
and all other minimal generators avoid $D$ and the last point.
Proposition~\ref{prop:universal-vertex-closure}, which allows the exceptional
outer core to strictly contain the active residual support, excludes
this final configuration.
\end{proof}

The required terminal-kernel estimate is as follows.
\begin{proposition}
\label{prop:universal-vertex-closure}
Let $t\ge2$, $a\ge3$, and $N\ge(t+1)(a-1)+2$.
Work on the ordered set $[N]$, put $U=[2,N]$, and let
$Y\subseteq X\subseteq U$ be initial intervals with
\[
 |X|=m\le a,\qquad |Y|=y=d+z+1,\qquad d\ge1,
 \quad z\ge0.
\]
Let $D_*$ consist of the first $d$ points of $Y$, and put
$Z=Y\setminus(D_*\cup\{\max Y\})$.
Let $\mathcal R$ be an upset kernel supported on $Z$, all of whose
minimal generators have size at least two. Its rank-$a$ extension
$\mathcal E_R$ is taken on $U\setminus D_*$ and is assumed shifted.
Define
\[
 \mathcal E=
 \{E\in\tbinom Ua:|E\cap Y|\ge2,\ E\cap D_*\ne\varnothing\}
 \;\dot\cup\;\mathcal E_R,
 \qquad
 \mathcal F=\{F\in\tbinom Ua:X\subseteq F\}.
\]
Let $\mathcal H,\mathcal G\subseteq\binom{[N]}a$ be the saturated
ordinary cross-intersecting pair specified by
\[
 \mathcal H(1)=\mathcal F^{\perp_{a-1}},\quad
 \mathcal H(\overline1)=\mathcal E,\qquad
 \mathcal G(1)=\mathcal E^{\perp_{a-1}},\quad
 \mathcal G(\overline1)=\mathcal F.
\]
Then
\[
 W_t(\mathcal H)W_t(\mathcal G)\le B^2,
\]
where
\[
 B=\max\{h,b\},\quad
 h=\binom{N-1}{a-1}-\binom{N-a-1}{a-1}+t,
 \quad b=(t+2)\binom{N-3}{a-2}+\binom{N-3}{a-3}.
\]
For $d\ge2$, the stronger bound $W_t(\mathcal H)W_t(\mathcal G)\le b^2$
holds.
\end{proposition}
\begin{proof}
Put $M=N-1$,
$S=\binom M{a-1}$, and write
\[
 P_m=S-\binom{M-m}{a-1},\quad
 f_m=\binom{M-m}{a-m},\quad
 J=\sum_{j=0}^{t-1}\binom t{j+1}c_{a+j}(\mathcal E),\quad
 Q=|\mathcal E^{\perp_{a-1}}|.
\]
The same clique count as in \eqref{eq:formula-s7} applies, since
every member of $\mathcal E$ meets $Y\subseteq X$ in at least two
points. Thus
\begin{equation}\label{eq:universal-routing-objective}
 W_t(\mathcal H)W_t(\mathcal G)=(P_m+J)(Q+t f_m).
\end{equation}

Both ordinary families have empty common intersection. For
$\mathcal H$, its section through $1$ has common intersection
$\{1\}$ before deleting $1$, since $|X|\ge2$ and $|U|\ge2a$;
its section avoiding $1$ is nonempty. For $\mathcal G$, every
$Y\setminus\{v\}$ with $v\in Y$ meets all generators of
$\mathcal E$ and extends to an $(a-1)$-set avoiding $v$.
The same extension can avoid any prescribed point outside $Y$.
Thus $\mathcal G(1)$ has empty common intersection, and
$\mathcal G(\overline1)$ is nonempty. The lifts are nontrivial by
Lemma~\ref{lem:clique-lift-duality}.

Suppose first that $d=1$ and $\mathcal R$ is empty. Then
$\{1,\min Y\}$ is a common transversal of $\mathcal H$ and
$\mathcal G$. Indeed, every member of $\mathcal H$ avoiding $1$
belongs to $\mathcal E$ and contains $\min Y$. Every member of
$\mathcal G$ avoiding $1$ contains $X$, and hence also $\min Y$. Let $C_*$ be the $(t-1)$-set used in the lift. An $(a+j)$-clique
in either ordinary family meets $\{1,\min Y\}$ in at least $j+1$
points, since otherwise deleting $j$ points gives an $a$-set
avoiding this transversal. Hence every lifted member meets
$C_*\cup\{1,\min Y\}$ in at least $t$ points, and
Theorem~\ref{thm:transversal} applies.
If $d=1$ and $\mathcal R$ is nonempty, then $z\ge2$ and
Proposition~\ref{prop:one-universal} applies. Its offset is
$c=\binom{M-y}{a-y}+t\binom{M-m}{a-m} \le(t+1)\binom{M-y}{a-y}$;
this satisfies the bound required in that proposition.

Now assume $d\ge2$. If $\mathcal R$ is nonempty, again $z\ge2$,
and Proposition~\ref{prop:multiple-universal} shows that the
product is less than $b^2$.

If $\mathcal R$ is empty and $z\ge1$, use
Proposition~\ref{prop:core-containment}. Every member of $\mathcal H$
meets $\{1\}\cup D_*$, so an $(a+j)$-clique has at least $j+1$
points there. Hence every lifted member has at least $t$ points in
$C_*\cup\{1\}\cup D_*$, a set of size $t+d$, where $C_*$ is the
$(t-1)$-set used in the lift. A set meeting every generating pair of $\mathcal E$ either contains
$D_*$ or omits exactly one point $v$ of $D_*$ and contains
$Y\setminus\{v\}$. Counting these disjoint possibilities for
$\mathcal G(1)$ and adding the $t f_m$ term gives
\[
 W_t(\mathcal G)
 =\binom{M-d}{a-d-1}+d\binom{M-y}{a-y}+t\binom{M-m}{a-m}
 \le\binom{M-d}{a-d-1}+(t+d)\binom{M-y}{a-y}.
\]
Here its parameter $w=y-d=z+1$ is at least two, so that proposition
shows that the product is less than $b^2$ for every $m\ge y$.

It remains to consider $d\ge2$, $\mathcal R=\varnothing$ and $z=0$.
Then $y=d+1\ge3$ and the graph on $Y$ is complete, since every pair
of $Y$ meets $D_*$. Thus $\mathcal E$ consists exactly of the $a$-sets
having at least two points of $Y$. Proposition~\ref{prop:unused-core-endpoints}
applied to \eqref{eq:universal-routing-objective} reduces the upper
bound to the two outer-support values $m=y$ and $m=a$.
For $m=y$, the ordinary families consist of the sets having at
least $2$ and at least $y$ points in $\{1\}\cup Y$, respectively.
An $(a+j)$-clique in such a threshold family raises its required
intersection by $j$. The lifted families thus require at least
$t+1$ and $t+y-1$ points, respectively, in the $(t+y)$-set
$C_*\cup\{1\}\cup Y$. This is the asymmetric threshold pair of
Lemma~\ref{lem:asymmetric-threshold}, so its product is strictly
below $b^2$.
If $y=a$, this also settles $m=a$.
For $m=a$ and $y<a$, the ordinary section $\mathcal F$ is a
singleton, and Lemma~\ref{lem:outer-clique-baseline} shows that
the product is less than $b^2$.
\end{proof}

\subsection{Core densities and universal vertices}\label{sec:kernel-estimates}
\label{sec:core-density}

We prove the estimates used in
Proposition~\ref{prop:universal-vertex-closure} by counting the
two-element generators and bounding the residual contribution by
the density inequality. Throughout, we use the following observation:
if $b\le|U|-a$, a $b$-set meets every member of an $a$-uniform family
on $U$ if and only if it meets every generator. Indeed, a disjoint
generator extends to an $a$-set in its complement.
\begin{lemma}\label{lem:min-weight}
Let $\mathcal L$ be a shifted $a$-uniform family on $[v]$, or on an
ordered ground set relabelled in increasing order. Then
\[
 c_{a+j}(\mathcal L)=\sum_{L\in\mathcal L}\binom{\min L-1}{j},
 \qquad
 W_t(\mathcal L)=\sum_{L\in\mathcal L}\binom{t+\min L-2}{t-1}.
\]
\end{lemma}
\begin{proof}
Since $\mathcal L$ is shifted, an $(a+j)$-set is a clique if and only
if its $a$ largest elements form a member of $\mathcal L$. For a fixed
member $L$ consisting of these largest elements, the remaining $j$
elements can be chosen in $\binom{\min L-1}{j}$ ways. This gives the
first formula. Multiplying by $\binom{t-1}{j}$ and summing over $j$,
we obtain the second formula by Vandermonde's identity.
\end{proof}

\begin{lemma}\label{lem:retained-core-density}
Let $a\ge3$, $t\ge2$, $d\ge0$, and
\[
 M\ge(t+1)(a-1)+1-d,\qquad
 r=a-d-1,\qquad q=a+s-1,
 \qquad1\le s\le\min\{t,r\}.
\]
If $\mathcal U\subseteq\binom{[M]}r$ and
$\mathcal V\subseteq\binom{[M]}q$ are cross $s$-intersecting with
densities $u$ and $v$, respectively, then
\[
 uv\le \sigma_s^2(1-u)(1-v),\qquad
 \sigma_s^2=
 \left(\frac{(M-s)_{s-1}}{(M)_{s-1}}\right)^2
 \frac{(r)_s(q)_s}{(M-r)_s(M-q)_s}.
\]
In particular, writing $h=|\mathcal U|$, $D=\binom Mr$, and
$\delta=1-h/D$, one has
\[
 h|\mathcal V|\le\delta C_s
       \binom{M-s}{r-s}\binom{M-s}{q-s},\qquad
 C_s=\frac{((M-s+1)_s)^2}{(M-r)_s(M-q)_s}.
\]
\end{lemma}
\begin{proof}
We use the operator from the proof of Lemma~\ref{lem:density}.
The singular-value formula \eqref{eq:density-singular-values} applies
whenever $r\le q$ and $M\ge r+q$. Here $r\le q$, and since $s\le t$,
we have $M-r-q\ge(s-1)(a-2)\ge0$.
The ratios of consecutive nonzero singular values are nonincreasing
with the degree, as in the proof of Lemma~\ref{lem:density}.
For $s=1$, the ratio is at most one since $M\ge r+q$.
For $s=r$, the only nonzero singular value on the mean-zero space is
the one in degree $s$. We may therefore assume $2\le s<r$.
After dropping a factor at most one, the squared ratio at degree $s$ is
at most
$\frac{s^2(a-d-1-s)(a-1)}{(M-a+d+1-s)(M-a-2s+1)}$.
The fraction decreases with $M$. For $M=(s+1)(a-1)+1-d$, the
difference between the denominator and numerator is
$(s-1)\bigl(as^2-2as-s^2+4s-1\bigr) +d(s-1)(sa-s+1)>0$.
For $s=2$ the first bracket equals $3$, and for $s\ge3$ it is positive
since $a\ge s+2$. Thus the norm on the mean-zero space is attained in degree $s$.
The centered-indicator argument in Lemma~\ref{lem:density} gives the
first inequality. Cancelling falling factorials as in
\eqref{eq:density-missing-factor} gives the second.
\end{proof}

\begin{proposition}\label{prop:one-universal}
Use the families and notation of
Proposition~\ref{prop:universal-vertex-closure}, with $d=1$ and
$y\ge3$. Then $W_t(\mathcal H)W_t(\mathcal G)\le B^2$,
with strict inequality if $\mathcal R\ne\varnothing$.
\end{proposition}

\begin{proof}
Put $x=2$, $K=N-2$, $A=a-1$, $L=N-1-y=K-z-1$, and
\[
 D=\binom K{a-2},\qquad
 c=\binom L{a-y}+t\binom{N-1-m}{a-m}
   \le(t+1)\binom L{a-y}.
\]
Let $\mathcal Q_R=\mathcal E_R^{\perp_{a-2}}$ on the $K$ points
$U\setminus\{x\}$ and put $h=|\mathcal Q_R|$.
A member of $\mathcal G(1)$ either contains $x$, with its remainder
in $\mathcal Q_R$, or omits $x$ and contains $Y\setminus\{x\}$.
The second condition is sufficient because $Y\setminus\{x\}$ meets
every residual generator. By the clique count leading to
\eqref{eq:formula-s7}, the only higher cliques of $\mathcal G$ are
$\{1\}\cup F$ for $F\in\mathcal F$. Thus $W_t(\mathcal G)=h+c$.
The family $\mathcal H$ is shifted: its section at $1$ consists of
the $(a-1)$-sets meeting the initial interval $X$, and inserting $1$
leaves a point of $Y\subseteq X$. In the residual section, inserting
$x$ leaves a point of its generator in $Z$, and shifts avoiding $x$
preserve $\mathcal E_R$. The sets containing $x$ and another point
of $Y$ are closed under shifts avoiding $1$.
Let $A_0$ be its weight when $\mathcal R=\varnothing$.
For a residual member of minimum $i$ in the increasing labeling of
$U\setminus\{x\}$, expand
$\binom{t+i}{t-1}=\sum_{j=0}^{t-1}\binom{t+1}{j+2}\binom{i-1}{j}$.
Lemma~\ref{lem:min-weight} then gives
\[
 W_t(\mathcal H)=A_0+J_R,\qquad
 J_R=\sum_{j=0}^{t-1}\binom{t+1}{j+2}c_{a+j}(\mathcal E_R).
\]

For the comparison pair with $\mathcal R=\varnothing$, we have
$h=D$ and $J_R=0$. The sections through $1$ are the companions of
the opposite sections avoiding $1$, and any two members avoiding
$1$ meet in $x$. Thus the pair is cross-intersecting. Each family
has a nonempty section avoiding $1$, and the two sections at $1$
have empty common intersection: the first comprises all
$(a-1)$-sets meeting $X$, and the second includes all $(a-1)$-sets
through $x$ as well as one containing $Y\setminus\{x\}$ and avoiding
$x$. Both families have common transversal $\{1,x\}$. Their lifts
with a common $(t-1)$-set $C$ are therefore nontrivial cross
$t$-intersecting families with common $t$-transversal
$C\cup\{1,x\}$, by the argument in the proof of
Proposition~\ref{prop:universal-vertex-closure}.
Theorem~\ref{thm:transversal} gives
\begin{equation}
 A_0(D+c)\le B^2. \label{eq:local-1-1}
\end{equation}
This settles the case $\mathcal R=\varnothing$. In the remaining
case, $z\ge2$ and $h<D$: a set disjoint from a fixed residual
member does not belong to $\mathcal Q_R$.
The families $\mathcal Q_R$ and $\mathcal K_{a+s-1}(\mathcal E_R)$
are cross $s$-intersecting: otherwise an $a$-subset of a clique
could avoid a member of $\mathcal Q_R$. Since $|Z|\le a-2$,
$\mathcal Q_R$ is nonempty. Thus $\mathcal K_{a+s-1}(\mathcal E_R)$
is empty for $s>a-2$. Applying Lemma~\ref{lem:density} with ranks
$a-2$ and $a+s-1$ on $K$ points and summing with the coefficients
in $J_R$, we obtain the following inequality, where $T$ is as in
Lemma~\ref{lem:scalar} with $M=K$:
\begin{equation}
 hJ_R\le(1-h/D)T. \label{eq:local-1-2}
\end{equation}
Lemma~\ref{lem:scalar} gives
\begin{equation}
 \frac T{b^2}<\frac t{t+1}(e-2),\qquad
 b=(t+2)\binom{N-3}{a-2}+\binom{N-3}{a-3}. \label{eq:local-1-3}
\end{equation}
For the first summand $T_1$ of $T$, the same lemma gives
\begin{equation}
 \frac{T_1}{b^2}\le\frac{t}{2(t+1)}. \label{eq:local-1-4}
\end{equation}

Every $(a-2)$-set containing at least $z-1$ points of $Z$ meets each
generator in $\mathcal R$, since each generator has size at least two.
Therefore
$h\ge z\binom{L+1}{a-y+1}+\binom{L+1}{a-y}$.
Since $K\ge(t+1)(a-1)$, cancellation gives
\begin{equation}
 \frac ch
 \le\frac{(t+1)(a-y+1)}{z(L+1)}
 =\frac{(t+1)(A-z)}{z(K-z)}\le\frac1z. \label{eq:local-1-5}
\end{equation}
Also $z\le A-1$ and
\begin{equation}
 \frac cD
 \le(t+1)\frac{(A-1)_z}{(K)_z}
                 \frac{K-A+1}{K-z}
 \le(t+1)\left(\frac{A-1}{K}\right)^z
 \le(t+1)^{1-z}. \label{eq:local-1-6}
\end{equation}

If $z\ge3$, equations \eqref{eq:local-1-3}, \eqref{eq:local-1-5}, and \eqref{eq:local-1-6} imply
\begin{align}
 \frac{(1+c/h)(1+c/D)T}{b^2}
 &<\left(1+\frac1z\right)
     \left(1+(t+1)^{1-z}\right)\frac t{t+1}(e-2)\\
 &\le\frac43\left(1+\frac1{(t+1)^2}\right)
                 \frac t{t+1}(e-2)
 <\frac43(e-2)<1. \label{eq:local-1-7}
\end{align}
If $z=2$, then $\mathcal R=\{Z\}$, and $\mathcal E_R$ is the complete
star through $Z$. This family has no clique of size greater than $a$,
so $T$ can be replaced by $T_1$ in \eqref{eq:local-1-2}. Equations \eqref{eq:local-1-4}--\eqref{eq:local-1-6} then give
\begin{equation}
 \frac{(1+c/h)(1+c/D)T_1}{b^2}
 \le\frac32\left(1+\frac1{t+1}\right)
              \frac{t}{2(t+1)}
 =\frac34\frac{t(t+2)}{(t+1)^2}<\frac34<1. \label{eq:local-1-8}
\end{equation}

Let $T_*=T$ for $z\ge3$ and $T_*=T_1$ for $z=2$.
By \eqref{eq:local-1-1}, \eqref{eq:local-1-2}, and either \eqref{eq:local-1-7} or \eqref{eq:local-1-8},
\begin{align*}
 W_t(\mathcal H)W_t(\mathcal G)
 &=A_0(h+c)+(h+c)J_R\\
 &\le\frac{h+c}{D+c}B^2
       +(1+c/h)\frac{D-h}{D}T_*\\
 &<\frac{h+c}{D+c}B^2
       +\frac{D-h}{D+c}B^2=B^2.
\end{align*}
\end{proof}

\begin{proposition}\label{prop:multiple-universal}
Use the families and notation of
Proposition~\ref{prop:universal-vertex-closure}, with $d\ge2$ and
$z\ge2$. Then $W_t(\mathcal H)W_t(\mathcal G)<b^2$, including
when $\mathcal R=\varnothing$.
\end{proposition}

\begin{proof}
Set $A=a-1$, $M=N-1$, $K=M-d$, $r=A-d$, $L=M-y=K-z-1$, and
\[
 S=\binom MA,\qquad D=\binom Kr,\qquad
 \mathcal Q_R=\mathcal E_R^{\perp_r},\quad h=|\mathcal Q_R|.
\]
Let
\begin{equation}\label{eq:local-2-11}
 c=d\binom L{r-z}+t\binom{M-m}{a-m}
   \le(t+d)\binom L{r-z}.
\end{equation}
A member of $\mathcal G(1)$ either contains $D_*$, with remainder
in $\mathcal Q_R$, or omits exactly one point of $D_*$ and contains
all other points of $Y$. Indeed, two omitted points of $Y$, one
in $D_*$, would form a disjoint generator. The second type contains
$Z$ and hence meets all residual generators. The higher cliques of
$\mathcal G$ are counted as in Proposition~\ref{prop:one-universal}.
The same shiftedness argument applies to $\mathcal H$: inserting
a point of $D_*$ in a residual member leaves a point of its
generator in $Z$, and shifts within the two-element-generator
part preserve both defining conditions. In the increasing
relabelling of $U\setminus D_*$, a residual member of minimum
$i$ has weight $\binom{t+d+i-1}{t-1}$; the coefficient of
$\binom{i-1}{s-1}$ in its expansion is $\binom{t+d}{d+s}$.
Thus Lemma~\ref{lem:min-weight} gives
\begin{equation}
 W_t(\mathcal G)=h+c,\qquad
 W_t(\mathcal H)=A_0+J_R,\qquad
 J_R=\sum_{s=1}^t\binom{t+d}{d+s}c_{a+s-1}(\mathcal E_R),                 \label{eq:local-2-1}
\end{equation}
where $A_0$ is the weight for $\mathcal R=\varnothing$.
When $\mathcal R=\varnothing$, $h=D$ and $J_R=0$.

Every $r$-set containing at least $z-1$ points of $Z$ meets each
generator in $\mathcal R$. Thus
$h\ge z\binom{L+1}{r-z+1}+\binom{L+1}{r-z}$.
Using $K\ge(t+1)A+1-d$, cancellation yields
\begin{equation}
 \frac ch\le\frac{(t+d)(r-z+1)}{z(K-z)}
      \le\frac{t+d}{(t+1)z},                                         \label{eq:local-2-2}
\end{equation}
because $(K-z)-(t+1)(r-z+1)\ge t(d+z-1)\ge0$.
Also $z\le r$ and
\begin{equation}
 \frac cD\le(t+d)\frac{(r)_z}{(K)_z}\frac{K-r}{K-z}
      \le\frac{t+d}{(t+1)^z}.                                       \label{eq:local-2-3}
\end{equation}
For the comparison with $\mathcal R=\varnothing$, the first
family $\mathcal H_0$ is generated by the pairs
$\{1,v\}$, $v\in X$, and the pairs in $Y$ meeting $D_*$.
Every such pair meets the initial $(d+1)$-set $\{1\}\cup D_*$.
Every $(a+j)$-clique contains at least $j+1$ points of
$\{1\}\cup D_*$. Hence each $k$-set in the clique lift contains at
least $t$ points of the fixed $(t+d)$-set formed by $\{1\}\cup D_*$
and the $t-1$ points used in the lift. Summing over the $t$-subsets
of this fixed set gives
\begin{equation}
 A_0\le\binom{t+d}d S.                                               \label{eq:local-2-4}
\end{equation}
Put $p=A/M$. Then
\[
 \frac DS=\frac{(A)_d}{(M)_d}\le p^d,\qquad
 \frac bS=p\left[t+2-(t+1)\frac{A-1}{M-1}\right]\ge(t+1)p,
 \qquad p\le\frac1{t+1}.
\]
Equations \eqref{eq:local-2-3}--\eqref{eq:local-2-4}, with $z\ge2$, give
\begin{equation}
 \frac{A_0(D+c)}{b^2}
 \le\frac{\binom{t+d}d}{(t+1)^d}
       \left(1+\frac{t+d}{(t+1)^2}\right)\le\frac{26}{27}<1.          \label{eq:local-2-5}
\end{equation}
The middle expression in \eqref{eq:local-2-5} decreases for $d\ge2$,
since the ratio of its value at $d+1$ to that at $d$ is at most
$(5/9)(14/13)<1$. At $d=2$, writing $u=1/(t+1)\le1/3$, the
expression is
$(1+u)(1+u+u^2)/2\le26/27$.

This proves the bound when $\mathcal R=\varnothing$.
For $\mathcal R\ne\varnothing$, the same bound on $A_0(D+c)$
will be used with the residual estimate. Here $0<h<D$, since
$r$-sets containing $Z$ are in $\mathcal Q_R$ whereas one disjoint
from a fixed residual member is not. The families $\mathcal Q_R$ and $\mathcal K_{a+s-1}(\mathcal E_R)$
are cross $s$-intersecting. Here the required intersection size is
$s$, whereas the lower index of the coefficient in
\eqref{eq:local-2-1} is $d+s$. For $s>r$, the clique family is empty. Lemma~\ref{lem:retained-core-density},
applicable since $K\ge(t+1)A+1-d$, gives
\begin{equation}
 hJ_R\le(1-h/D)T_d,\qquad
 T_d=\sum_{s=1}^{\min(t,r)}\binom{t+d}{d+s}E_s(d),                    \label{eq:local-2-6}
\end{equation}
where, for $K_d=M-d$ and $r_d=A-d$,
\[
 E_s(d)=\binom{K_d}{r_d}\binom{K_d}{A+s}
 \left(\frac{(K_d-s)_{s-1}}{(K_d)_{s-1}}\right)^2
 \frac{(r_d)_s(A+s)_s}{(K_d-r_d)_s(K_d-A-s)_s}.
\]
Let $C=M-1-A\ge tA$ and $x=A/C\le1/t$.
Direct cancellation gives, whenever both sides exist,
\begin{equation}
 \frac{E_s(d+1)}{E_s(d)}=
 \frac{(A-d-s)(K_d-A-s)^2(K_d-2s+1)^2}
 {(K_d-A-2s)(K_d-s)^2(K_d-s+1)^2}
 \le\frac{A-d-s}{K_d-A-2s}\le x.                                   \label{eq:local-2-7}
\end{equation}
All denominators are positive. The last inequality follows from
$C(d+s)\ge A(d+2s-1)$, since $C\ge2A$.
The calculation in Lemma~\ref{lem:scalar}, with its parameters
$(A,M)=(a-1,N-2)$, gives
\begin{equation}
 E_s(1)\le x^{s-1}E_1(1),\qquad
 \frac{E_1(1)}{b^2}
 =\frac{A-1}{A}\frac C{C-1}
   \left[t+2-(t+1)\frac{A-1}{M-1}\right]^{-2}
 \le\frac1{(t+1)^2}.                                                \label{eq:local-2-8}
\end{equation}
Equations \eqref{eq:local-2-7}--\eqref{eq:local-2-8} imply
\begin{equation}
 \frac{T_d}{b^2}\le
 f_{t,d}:=\sum_{s=1}^t
 \frac{\binom{t+d}{d+s}}{(t+1)^2t^{d+s-2}}
 <\frac{(d+2)^2}{9(d+1)2^{d-1}}.                                   \label{eq:local-2-9}
\end{equation}
Indeed the first summand is
\[
 \frac1{(d+1)!}\frac t{t+1}\prod_{i=2}^d\left(1+\frac it\right)
 \le\frac{d+2}{9\,2^{d-1}}.
\]
For $d\ge2$, the factors $t/(t+1)$ and $1+2/t$ have product
$(t+2)/(t+1)$, which increases with $1/t$. The remaining factors
also increase with $1/t$, so the expression is largest at $t=2$.
Since the ratio of consecutive summands is
$(t-s)/(t(d+s+1))<1/(d+2)$, a geometric sum gives
\eqref{eq:local-2-9}.
From \eqref{eq:local-2-2}--\eqref{eq:local-2-3}, $z\ge2$ and $t\ge2$,
\[
 1+c/h\le1+\frac{d+2}{6},\qquad 1+c/D\le1+\frac{d+2}{9}.
\]
Consequently
\begin{equation}
 \frac{(1+c/h)(1+c/D)T_d}{b^2}
 <\frac{(d+8)(d+11)(d+2)^2}{486(d+1)2^{d-1}}
 \le\frac{520}{729}<1.                                              \label{eq:local-2-10}
\end{equation}
The rational expression in \eqref{eq:local-2-10} decreases for
$d\ge2$. Indeed, the ratio of its value at $d+1$ to its value at $d$
is at most
$\tfrac12\cdot\tfrac{11}{10}\cdot\tfrac{14}{13}\cdot(\tfrac54)^2
=385/416<1$, even after dropping the factor $(d+1)/(d+2)$.

Combining \eqref{eq:local-2-5}, \eqref{eq:local-2-6}, and \eqref{eq:local-2-10},
\[
 (A_0+J_R)(h+c)
 \le\frac{h+c}{D+c}b^2+(1+c/h)\frac{D-h}{D}T_d
 <\frac{h+c}{D+c}b^2+\frac{D-h}{D+c}b^2=b^2.
\]
\end{proof}

\subsection{A core-containment estimate}\label{sec:core-containment}
\begin{proposition}\label{prop:core-containment}
Let $t\ge2$, $k\ge t+2$, $n\ge(t+1)(k-t+1)$, and
$2\le d\le k-t$. Put $q=t+d$, fix a $q$-set $Q$, and let
$2\le w\le n-q$.
Suppose
\[
 \mathcal H\subseteq\{H\in\tbinom{[n]}k:|H\cap Q|\ge t\},
\]
and
\[
 |\mathcal G|\le D_q+qE_w,\qquad
 D_q=\binom{n-q}{k-q},\quad
 E_w=\binom{n-q-w}{k-q-w+1}.
\]
Then $|\mathcal H||\mathcal G|<b_t(n,k)^2$.
\end{proposition}

\begin{proof}
Put $M=n-t$, $A=k-t$, $S=\binom MA$, $p=A/M$ and $z=p(1-p)$.
Thus $p\le1/(t+1)\le1/3$ and $M\ge3A$.
The sequence $E_w$ is nonincreasing in $w$, since each step
decreases the two binomial indices by one. Thus
\begin{equation}
 \frac{E_w}{D_q}\le\frac{E_2}{D_q}
 =\frac{(A-d)(M-A)}{(M-d)(M-d-1)}\le p(1-p)=z. \label{eq:local-3-1}
\end{equation}
For the last comparison it suffices to check
$(A-d)M^2\le A(M-d)(M-d-1)$.
The right side minus the left equals
$dM^2-A(2d+1)M+Ad(d+1) \ge AM(d-1)+Ad(d+1)>0$,
using $M\ge3A$ and $d\ge2$.
If $A=d$, then $E_2=0$ and \eqref{eq:local-3-1} is immediate.

Counting the members of $\mathcal H$ according to $|H\cap Q|$ gives
\begin{align*}
 \frac{D_q|\mathcal H|}{S^2}
 &\le\sum_{j=0}^d\binom{t+d}{t+j}
       \frac{(A)_d(A)_j(M-A)_{d-j}}{(M)_d^2}\\
 &\le\sum_{j=0}^d\binom{t+d}{t+j}
       p^{d+j}(1-p)^{d-j}=:A_{d+1}.
\end{align*}
Indeed, for $0\le v<d$, the inequality
$1-v/A\le(1-v/M)^2$ follows from $A/M\le1/2$.
Multiplication gives $(A)_d/A^d\le((M)_d/M^d)^2$; the other
normalized numerator factors are at most one.
Equivalently, for every $r\ge3$,
\[
 A_r=\sum_{i=0}^{r-1}\binom{t+r-1}i p^{2(r-1-i)}z^i.
\]
Combining with \eqref{eq:local-3-1},
\begin{equation}
 \frac{|\mathcal H||\mathcal G|}{S^2}
 \le(1+(t+d)z)A_{d+1}. \label{eq:local-3-2}
\end{equation}
We have
$A_{r+1}=pA_r+\binom{t+r-1}r z^r \le\left(p+\frac{tz}{r}\right)A_r$.
For the last inequality, use
$A_r\ge\binom{t+r-1}{r-1}z^{r-1}$.
Since $p\le1/3$, $tz<1$ and $z\le1/4$, for $d\ge2$ the ratio of
the right-hand side of \eqref{eq:local-3-2} at $d+1$ to its value at
$d$ is at most
$(1+z)\left(p+\frac{tz}{d+1}\right) <\frac54\left(\frac13+\frac13\right)=\frac56$.
Thus $d=2$ gives the largest upper bound.
At $d=2$, division by $(t+2)^2z^2$ gives
$(1+(t+2)z) \left[\frac{p^2}{(t+2)^2(1-p)^2} +\frac{p}{(t+2)(1-p)}+\frac{t+1}{2(t+2)}\right]$.
Every factor is increasing for $0<p\le1/(t+1)$.
At $p=1/(t+1)$, one minus this expression equals
$\frac{2t^4+3t^3-6t^2-10t-2} {2t^2(t+1)(t+2)^2}>0$.
Indeed, putting $u=t-2\ge0$, we can write the numerator as
$2u^4+19u^3+60u^2+66u+10$.
Finally
$\frac{b_t(n,k)}S =(t+2)\frac{\binom{M-2}{A-1}}{\binom MA} +\frac{\binom{M-2}{A-2}}{\binom MA} \ge(t+2)z$.
The right-hand side of \eqref{eq:local-3-2} is therefore strictly
smaller than $b_t(n,k)^2/S^2$.

\end{proof}

\subsection{Unused support points}\label{sec:unused-support}
\begin{proposition}\label{prop:unused-core-endpoints}
Let $t\ge2$, $a\ge3$, $M\ge(t+1)(a-1)+1$, and
$3\le y\le a$. Let $\mathcal E\subseteq\binom{[M]}a$ be generated
by subsets of $Y=[y]$, each of size at least two. Put
\[
 J=\sum_{j=0}^{t-1}\binom t{j+1}c_{a+j}(\mathcal E),\qquad
 Q=|\mathcal E^{\perp_{a-1}}|,\qquad S=\binom M{a-1}.
\]
For $y\le m\le a$, set
\[
 P_m=S-\binom{M-m}{a-1},\quad
 f_m=\binom{M-m}{a-m},\quad
 V_m=(P_m+J)(Q+t f_m).
\]
Then
\[
 V_m\le\max\{V_y,V_a\}\qquad(y\le m\le a).
\]
Thus the maximum occurs when the exceptional core is $Y$, or when
it has size $a$ and the exceptional family is a singleton.
\end{proposition}

\begin{proof}
Every $(a-1)$-set containing at least $y-1$ points of $Y$ meets every
generator, since generators have size at least two. Consequently
\begin{equation}
 Q\ge y\binom{M-y}{a-y}\ge3f_m\qquad(m\ge y), \label{eq:local-4-1}
\end{equation}
since $f_m$ decreases with $m$.
For $y\le m<a$, define
\[
 \delta_m=P_{m+1}-P_m=\binom{M-m-1}{a-2},\qquad
 \gamma_m=t(f_m-f_{m+1})=t\binom{M-m-1}{a-m}.
\]
and put $R_m=\delta_m/\gamma_m$.
Then
\[
 V_{m+1}-V_m
 =\gamma_m g_m,\qquad
 g_m=R_m Q-J+C_m,\qquad
 C_m=tR_m f_m-P_{m+1}.
\]
For $m\le a-2$, direct binomial cancellation gives
\[
 \frac{R_{m+1}}{R_m}
 =\frac{M-a-m+1}{a-m},\qquad
 C_{m+1}-C_m
 =-\delta_m\frac{a-1}{M-a}-\delta_{m+1}.
\]
Hence $g_{m+1}\ge g_m$ whenever
\begin{equation}
 Q\ge\gamma_m\frac{a-m}{M-2a+1}
       \left(\frac{a-1}{M-a}+\frac{\delta_{m+1}}{\delta_m}\right). \label{eq:local-4-2}
\end{equation}
The right side is strictly less than $3f_m$: indeed,
$\gamma_m\le t f_m,\quad M-2a+1\ge(t-1)(a-1),\quad \frac{a-1}{M-a}\le\frac1t,\quad \frac{\delta_{m+1}}{\delta_m}\le1$,
and $a-m<a-1$. It is therefore at most
$\frac{t+1}{t-1}\frac{a-m}{a-1}f_m<3f_m$.
Equation \eqref{eq:local-4-1} proves \eqref{eq:local-4-2}, strictly.
Hence $g_m$ increases with $m$, so a positive difference
$V_{m+1}-V_m$ cannot be followed by a negative one. The maximum of
$V_m$ is therefore attained at an endpoint.
\end{proof}

\subsection{Threshold families}\label{sec:explicit-endpoints}
\begin{lemma}\label{lem:asymmetric-threshold}
Let $t\ge2$, $k\ge t+2$ and $n\ge(t+1)(k-t+1)$.
For $d\ge2$ with $t+d\le k$, put $L=t+d+1$ and
\[
 \mathcal A_d=\{F\in\tbinom{[n]}k:|F\cap[L]|\ge t+1\},\qquad
 \mathcal B_d=\{F\in\tbinom{[n]}k:|F\cap[L]|\ge t+d\}.
\]
Then $|\mathcal A_d||\mathcal B_d|<b_t(n,k)^2$.
\end{lemma}
\begin{proof}
Write
\[
 K(L,q)=\sum_{j\ge q}\binom Lj\binom{n-L}{k-j},\qquad
 F_r=K(t+2r,t+r).
\]
For fixed $L$, the sequence $K(L,q)$ is log-concave in $q$.
Indeed, its summands $x_j=\binom Lj\binom{n-L}{k-j}$ form a
log-concave sequence whose positive terms have consecutive indices.
Writing $T_q=\sum_{j\ge q}x_j$, we have
$T_q^2-T_{q-1}T_{q+1}=x_qT_q-x_{q-1}T_{q+1}\ge0$,
since $x_qx_{q+i-1}\ge x_{q-1}x_{q+i}$ for every $i\ge1$.
Thus, for a fixed sum of the thresholds, the product of two tails
does not decrease when the thresholds are moved closer together.

We also use the following inequality for blocks of different sizes,
with $1\le q<k$ and $q<L\le n-k+q$:
\begin{equation}\label{eq:mixed-block-tails}
 K(L,q)K(L,q+1)\le K(L-1,q)K(L+1,q+1).
\end{equation}
Put $A=K(L-1,q)$ and $V=K(L,q+1)$. Counting the sets that
contain the newly added point gives
$K(L,q)=A+\alpha,\quad K(L+1,q+1)=V+\beta$,
where
$\alpha=\binom{L-1}{q-1}\binom{n-L}{k-q},\quad \beta=\binom Lq\binom{n-L-1}{k-q-1}$.
Both quantities are positive in the stated range. For $j\ge q$,
compare the summand of $V$ at index $j+1$ with that of $A$ at index
$j$, whenever the former is nonzero:
\[
 \frac{\binom L{j+1}\binom{n-L}{k-j-1}}
      {\binom{L-1}j\binom{n-L+1}{k-j}}
 =\frac{L(k-j)}{(j+1)(n-L+1)}
 \le\frac{L(k-q)}{(q+1)(n-L+1)}
 <\frac{L(k-q)}{q(n-L)}=\frac\beta\alpha.
\]
Thus $\alpha V\le\beta A$, proving
\eqref{eq:mixed-block-tails}.

If $d$ is odd, write $d=2r-1$. Tail log-concavity bounds the desired
product by $F_r^2$. If $d=2r$ is even, it bounds the product by
$K(t+2r+1,t+r)K(t+2r+1,t+r+1)$; equation
\eqref{eq:mixed-block-tails} bounds this by $F_rF_{r+1}$.
Indeed, here $L=t+2r+1$ and $q=t+r$, and $2r\le k-t$ gives
$q<k$ and $n-L\ge k-q$.

It remains to compare the $F_r$. We have $F_1=b_t(n,k)$.
For $1\le r\le k-t$, put $Q=t+2r$. In this range $Q+2\le n$.
Counting according to the two points added to the block gives
$F_r-F_{r+1} =\binom Q{t+r}\binom{n-Q-2}{k-t-r} -\binom Q{t+r-1}\binom{n-Q-2}{k-t-r-1}$.
If $r=k-t$, the second term is zero. Otherwise, cancellation
shows that the first term is at least the second provided
$n\ge\left(2+\frac{t-1}{r+1}\right)(k-t+1)$,
which follows from the hypothesis on $n$. Thus $F_r$ is
nonincreasing for $1\le r\le k-t+1$.
Moreover, $F_1>F_2$: for $r=1$ this threshold is strictly
smaller than $(t+1)(k-t+1)$, and the first binomial product is positive.
If $d$ is odd, then $r\ge2$ and $F_r^2<F_1^2$. If $d$ is even,
then $r+1\ge2$ and $F_rF_{r+1}<F_1^2$. The bound also holds
when one of the families is empty.
\end{proof}

\medskip\label{sec:outer-clique}
The remaining endpoint has a complete graph as its kernel.
\begin{lemma}
\label{lem:outer-clique-baseline}
Let $t\ge2$, $a\ge3$, $N\ge(t+1)(a-1)+2$, and put $M=N-1$.
Let $Y$ be an initial $(d+1)$-set in an $M$-point set, where
$2\le d\le a-1$, and set
\[
 \mathcal E=\{E\in\tbinom{[M]}a:|E\cap Y|\ge2\},\quad
 V=|\mathcal E^{\perp_{a-1}}|,\quad
 J=\sum_{j=0}^{t-1}\binom t{j+1}c_{a+j}(\mathcal E).
\]
Write
\[
 S=\binom M{a-1},\quad
 u_a=S-\binom{M-a}{a-1},\quad
 b=(t+2)\binom{M-2}{a-2}+\binom{M-2}{a-3}.
\]
Then $(u_a+J)(V+t)<b^2$. If $d\le a-2$, the same strict bound
holds with $S$ in place of $u_a$.
\end{lemma}
\begin{proof}
If $d=a-1$, take $X=Y$ in the section formulas of
Proposition~\ref{prop:universal-vertex-closure}, so $f_a=1$.
The count in \eqref{eq:universal-routing-objective} is valid for
these sections regardless of saturation and gives the product
$(u_a+J)(V+t)$. The two lifts are precisely the threshold families
in Lemma~\ref{lem:asymmetric-threshold}, as in the $m=y$ case of
Proposition~\ref{prop:universal-vertex-closure}. The lemma gives the
required strict bound.
Assume henceforth $d\le a-2$ and put
\[
 p=\frac{a-1}{M}\le\frac1{t+1},\qquad
 D=\binom{M-d}{a-d-1},\qquad
 e=\binom{M-d-1}{a-d-1}.
\]
A set meeting every pair of $Y$ contains at least $d$ points of $Y$.
Consequently $V=D+de$. Moreover,
\[
 \frac DS=\frac{(a-1)_d}{(M)_d}\le p^d,
 \qquad
 D\ge M-d\ge t(d+1)+2.
\]
In particular, with $g_d=d+1+1/(d+1)$,
\begin{equation}\label{eq:outer-clique-G}
 V+t<(d+1+1/(d+1))D\le g_dp^dS.
\end{equation}

Every $(a+j)$-clique of $\mathcal E$ contains at least $j+2$ points
of $Y$. A union bound therefore gives
$c_{a+j}(\mathcal E)\le\binom{d+1}{j+2}\binom{M-j-2}{a-2} \quad(0\le j\le d-1)$,
and there are no larger cliques. For $0\le j\le d-1$,
\[
 \frac{\binom{M-j-2}{a-2}}S
 =\frac{a-1}{M-j-1}
       \prod_{i=0}^j\frac{M-a+1-i}{M-i}
 \le p(1-p)^j.
\]
Thus
\begin{equation}\label{eq:outer-clique-J}
 \frac JS\le p\sum_{r=1}^{\min\{t,d\}}
        \binom tr\binom{d+1}{r+1}(1-p)^{r-1}.
\end{equation}
We also have
$\frac bS =p\left(t+2-\frac{(t+1)(a-2)}{M-1}\right) \ge(t+1)p$.
Combining these estimates bounds $(S+J)(V+t)/b^2$ by
\[
 \frac{g_d}{(t+1)^2}
 \left[p^{d-2}+\sum_{r=1}^{\min\{t,d\}}
  \binom tr\binom{d+1}{r+1}p^{d-1}(1-p)^{r-1}\right].
\]
Each term in brackets is nondecreasing for $0\le p\le1/(t+1)$.
Indeed, $p^{d-1}(1-p)^{r-1}$ is nondecreasing up to
$p=(d-1)/(d+r-2)\ge1/2$. Substituting $p=1/(t+1)$, replacing
$(1-p)^{r-1}$ by $1$, and summing, we obtain the upper bound
\begin{equation}\label{eq:outer-clique-Q}
 Q_d(t)=\frac{g_d}{(t+1)^{d+1}}
       \left[\binom{t+d+1}d+t-d\right].
\end{equation}

We show that $Q_d(t)$ decreases with $d$ when $t\ge3$ and
$d\ge2$, or when $t=2$ and $d\ge3$. Put
$B_d=\binom{t+d+1}d+t-d$. The elementary bound
$\binom{t+d+1}d\ge(t+1)(d+1)$, which holds at $d=2$ and follows
for larger $d$ by induction, implies
\[
 B_{d+1}\le\frac{t+d+3}{d+1}B_d,
 \qquad \frac{g_{d+1}}{g_d}\le\frac{d+2}{d+1}.
\]
For the first inequality, multiply by $d+1$ and subtract the
left-hand side from the right-hand side to obtain
$\binom{t+d+1}d+(t+2)(t-d)+d+1\ge(t+2)(t+1)$.
Therefore
$\frac{Q_{d+1}(t)}{Q_d(t)} \le\frac{(d+2)(t+d+3)}{(d+1)^2(t+1)}<1$
in the stated ranges, since the difference between the denominator
and numerator is $t(d^2+d-1)-3d-5>0$.
For $t\ge3$,
$Q_2(t)=\frac{5(t^2+7t+2)}{3(t+1)^3}\le Q_2(3)=\frac56$,
where the inequality follows by differentiation: the sign of the
derivative is the sign of $-t^2-12t+1$.
For $t=2,d\ge3$, $Q_d(2)\le Q_3(2)=323/324$.
This proves all cases except $t=d=2$.

For $t=d=2$, put $A=a-1\ge3$ and write
\[
 X=\binom{M-3}{a-2},\qquad Y_0=\binom{M-3}{a-3}.
\]
Counting the cliques gives $J=7X+2Y_0$. We have
\[
 \frac XS\le p(1-p),\qquad \frac{Y_0}S\le p^2.
\]
Cancelling binomial coefficients reduces the first inequality to
$A\ge2-2/M$. For the second, $Y_0/S$ is the proportion of $A$-sets
that contain two fixed points and avoid a third. This is at most
$(A)_2/(M)_2\le p^2$.
If $a\ge5$, then $D\ge M-2\ge10$, and
$V+2\le3D+2\le16D/5$. If $a=4$, direct counting gives
$V+2=3D$, so the same bound holds. Since $D/S\le p^2$,
$\frac{(S+J)(V+2)}{b^2} \le\frac{16}{5}\, \frac{1+7p-5p^2}{(4-3p)^2} \le\frac{80}{81}<1$.
The last inequality follows since the derivative of the rational
function has positive numerator $34-19p$ for $0\le p\le1/3$.
\end{proof}

\section{Equality cases and proof of the main theorem}\label{sec:equality}
Throughout this section, $t\ge2$, $k\ge t+2$ and
$n\ge(t+1)(k-t+1)$. In particular, writing $r=k-t$ gives
$n-2k\ge(t-1)(r-1)>0$. The terminal-pair reduction,
Theorem~\ref{thm:critical} and Proposition~\ref{prop:shifted-bound}
give the upper bound $m_t(n,k)^2$; every maximizing pair is saturated.
We classify pairs with a common $t$-transversal of size $t+1$, prove
that every shifted extremal pair has one, and reverse the shifts.

\subsection{Extremal pairs after shifting}
\begin{lemma}\label{lem:transversal-equality}
Let $\mathcal A,\mathcal B$ be nontrivial cross $t$-intersecting
families with a common $t$-transversal of size $t+1$. If
$|\mathcal A||\mathcal B|=m_t(n,k)^2$, then the pair is of type
\textup{(i)} or \textup{(ii)} in Theorem~\ref{thm:main}.
\end{lemma}
\begin{proof}
Relabel the transversal as $S=[q]$ and retain $q,r,N,a,b,D,M$
and $g$ from the proof of Theorem~\ref{thm:transversal}.
By saturation, both families
contain the full $S$-star. Let $x\ge y$ be the two largest section
sizes of $\mathcal A$. Equality must hold in each comparison following
\eqref{eq:transversal-product}.

The equality statement in Lemma~\ref{lem:offset-star-companion} gives
\[
 (a+qs)(a+qg(s))<M^2\qquad(1\le s<b).
\]
At $s=b$, equality requires $M=a+qb$.
The comparison of $P(x,y)$ with its two diagonal endpoints is also
strict when $x>y$, since
$(q-2)a-x+(q-1)^2y>0$ for $a>b\ge x$ and $q\ge3$.
Thus if $x,y\le b$, then $x=y=b$. If $x,y\ge b$, the same
conclusion follows by considering $g(y)$ and $g(x)$. Here $g(s)<b$
for $s>b$: the companion is contained in the full star and an
additional set avoiding the centre excludes some members of that star.

In the remaining case $x\ge b\ge y$, put $X=x$ and $Y=g(y)$ as
in the proof of Theorem~\ref{thm:transversal}. Equality requires
$g(Y)=y$ and
\begin{equation}\label{eq:equality-single-family-sizes}
 a+X+(q-1)g(X)=a+Y+(q-1)g(Y)=M.
\end{equation}
For $z\in\{X,Y\}$, the corresponding nontrivial $t$-intersecting
family attains the bound in Theorem~\ref{thm:nontrivial-ak}. Its
section sizes with respect to $S$ are $z,g(z),\ldots,g(z)$.

For a radius-one family, every $t$-transversal of size $t+1$
is contained in its defining $(t+2)$-set, and all section sizes are $b$.
For $\mathcal H(T,Z)$, every such transversal is $T\cup\{u\}$ with
$u\in Z$, and the section sizes are $D,1,\ldots,1$. Indeed, since
$n>2k$, a transversal of size $t+1$ must meet every complete-star
generator in at least $t$ points. The radius-one generators are all
the $(t+1)$-subsets of the defining block. The Hilton--Milner generators include
$T\cup\{z\}$ for every $z\in Z$. A $(t+1)$-set omitting at least
two points of $T$ cannot meet these
generators in $t$ points. If it omits exactly one, it must contain
all of $Z$, requiring at least $(t-1)+(r+1)>t+1$ points. Thus it
contains $T$, and the exceptional members force its remaining point
into $Z$.
It follows from~\eqref{eq:equality-single-family-sizes} that
$X,Y\in\{b,D\}$, where $g(b)=b$ and $g(D)=1$ in the
Hilton--Milner case. If $X\ne Y$, then
$(q-1)(X-Y)(g(X)-g(Y))<0$, contradicting equality in the
rearrangement comparison. After permuting the omitted labels,
the only possible section sizes are
\begin{equation}\label{eq:equality-section-patterns}
 (b,b,\ldots,b)\quad\hbox{on both sides},\qquad
 (D,1,\ldots,1)\quad\hbox{on both sides}.
\end{equation}
Equality in~\eqref{eq:transversal-product} forces every section
to attain its bound. In the second pattern, $x=D$ and $y=1$,
so on the opposite side only the section with the first omitted
label can have size $g(1)=D$; all others have size $g(D)=1$.
Thus the large sections have the same omitted label.

We recover the sections before lexicographic replacement.
In the first pattern, consider two sections on opposite sides with
different omitted labels. They are cross-intersecting and both have
density $r/N$. For these densities equality holds in
\eqref{eq:transversal-density}. In its proof, the singular values of
the disjointness operator strictly decrease with the level. Equality
in the norm bound therefore puts the centred indicators of both
sections in the first level. Each indicator is thus of the form
$c+\sum_{i\in R}\lambda_i$ on the $r$-sets $R$ of the outside ground
set.

Comparing two $r$-sets that differ by one point shows that
$\lambda_i-\lambda_j\in\{-1,0,1\}$ for every $i,j$. The coefficients
take two distinct values, since the indicator is nonconstant.
If each value occurred at least twice, two $r$-sets with the same
$r-2$ other points would give indicator values differing by two.
Hence one coefficient value occurs just once, so the section is a
full star or the complement of one. Its density $r/N<1/2$ forces
a full star. The two stars must have the same centre, since
stars with distinct centres contain disjoint members when $N\ge2r$.

The bipartite graph on the omitted labels, with edges between
unequal labels, is connected because $q\ge3$. Hence every section
on both sides has the same centre $u$, and both original families are
$\mathcal A_1(S\cup\{u\})$.

In the second pattern, let $v\in S$ be the common omitted label of
the large sections and put $T=S\setminus\{v\}$. Write the singleton
sections of the first side as $\{R_i\}$ and those of the second as $\{Q_i\}$,
$i\in T$, where all these sets have size $r$ and avoid $S$.
The large first section is contained in the full companion of each
$Q_i$, whose size is exactly $D$. Hence these companions are equal,
and all $Q_i$ are equal to some $r$-set $Q$. Indeed, distinct
$r$-sets have distinct full companions when $N>2r$: a point in one
set but not the other extends to an $r$-set disjoint from the latter.
Similarly, all $R_i$ are equal to some $r$-set $R$. The large first
section is the full companion of $Q$, and the large second section
is the full companion of $R$.
As $|T|=t\ge2$, two different omitted labels give
$R\cap Q\ne\varnothing$. Taking
$X=R\cup\{v\}$ and $Y=Q\cup\{v\}$ recovers precisely
$(\mathcal H(T;X,Y),\mathcal H(T;Y,X))$, with $|X\cap Y|\ge2$.
The sizes of the recovered families give the required comparisons
between $h_t(n,k)$ and $b_t(n,k)$.
\end{proof}

We now show that shifted extremal pairs have the common transversal
required above.
\begin{lemma}\label{lem:strict-kernel-equality}
Under the hypotheses of Proposition~\ref{prop:universal-vertex-closure},
suppose that the blocker of $\mathcal E$ has no singleton generator.
Then $W_t(\mathcal H)W_t(\mathcal G)<B^2$.
\end{lemma}
\begin{proof}
The case $d=1$, $\mathcal R=\varnothing$ is impossible: every member
of $\mathcal E$ then contains the unique universal vertex, giving a
singleton generator of its blocker. For $d=1$ with
$\mathcal R\ne\varnothing$, Proposition~\ref{prop:one-universal}
is strict, also when $m>y$.
For $d\ge2$, all cases in the proof of
Proposition~\ref{prop:universal-vertex-closure} are strict by
Propositions~\ref{prop:multiple-universal} and
\ref{prop:core-containment}, except possibly the complete graph on
the inner support. For that case, the endpoint with outer support
equal to the inner support is strict by
Lemma~\ref{lem:asymmetric-threshold}; the other endpoint, if distinct,
is strict by Lemma~\ref{lem:outer-clique-baseline}.
Proposition~\ref{prop:unused-core-endpoints} completes the proof.
\end{proof}

\begin{lemma}\label{lem:shifted-equality-cover}
Let $\mathcal A$ and $\mathcal B$ be nonempty nontrivial cross
$t$-intersecting families. If both families are shifted and
$|\mathcal A||\mathcal B|=m_t(n,k)^2$, then they have a common
$t$-transversal of size $t+1$.
\end{lemma}
\begin{proof}
Suppose otherwise. The pair is extremal and saturated.
If each family has a canonical generator of size $t+1$, then,
by shiftedness, each contains the full $[t+1]$-star. The generators
on opposite sides meet in at least $t$ points, so $[t+1]$ is a
common $t$-transversal. We may therefore assume that every canonical
generator of $\mathcal B$ has size at least $t+2$.

We follow the terminal exchanges in the proof of
Proposition~\ref{prop:shifted-bound}. The paired ranks exclude
$s=t+2$ under our assumption on the generators of $\mathcal B$. The case $s=t+3$
is excluded by the strict bound in Lemma~\ref{lem:asymmetric-threshold}. For $s\ge t+4$, both exchanges
at a pair of middle ranks preserve nontriviality, by
Lemma~\ref{lem:terminal-deletion}. The numerical lemmas in
Appendix~\ref{sec:exchange-numerics} then give a strict increase
in one direction, with two exceptions when $t=2$. In those cases
either Lemma~\ref{lem:t2-global-middle} gives a strict increase or
the pair consists of two complete-intersection families of radius
two or three. The latter pairs have product strictly less than
$b_2(n,k)^2$, since $F_1>F_2\ge F_3$ in the proof of
Lemma~\ref{lem:asymmetric-threshold}. The remaining active ranks
are $t+1$ and $s-1$. The retained B-family has a common $t$-set,
since it is cross $t$-intersecting with the full star through a
shortened active generator of $\mathcal A$. This gives \eqref{eq:formula-s4}
and \eqref{eq:formula-s6} as before.

In those sections, the blocker $\mathcal P$ of $\mathcal E$ has
no singleton generator, since $\mathcal B$ has no canonical
generator of size $t+1$. If $c>t$, Lemma~\ref{lem:strict-kernel-equality}
therefore excludes the resulting universal-vertex kernel.
If $c=t$, retain the actual inner support $Y$ and deform the outer
core as in~\eqref{eq:formula-s7}. The deformation preserves
$\mathcal E$ and $\mathcal P$, including the absence of singleton
generators. The endpoint inequality~\eqref{eq:formula-s8}, together
with the global upper bound, makes at least one endpoint extremal.
At the endpoint with outer core $Y$, the B-family still has
no generator of size $t+1$, as in the proof of
Proposition~\ref{prop:shifted-bound}. If $|Y|=3$, the strict
asymmetric threshold bound applies. If $|Y|\ge4$, the
terminal-generator argument applies to this extremal pair.
The last point of $Y$ is in a minimal generator of $\mathcal E$,
so the new terminal configuration has $c>t$. It gives a kernel
with a universal vertex, whose product is strictly below the
bound by Lemma~\ref{lem:strict-kernel-equality}. At the singleton endpoint
with outer core of size $a$, the inner exchanges
\eqref{eq:formula-s9}--\eqref{eq:formula-s11}
strictly increase the product unless they give such a kernel.
The exact cylinder counts and preservation of nontriviality apply
at equality as well as above the bound. The resulting kernel is strict
by the same lemma, excluding the second endpoint too.
\end{proof}

\subsection{Recovering a pair from one shift}
We use a bipartite graph to describe the possible inverse shifts. Write $P+i$ for $P\cup\{i\}$ when $i\notin P$.

\begin{lemma}\label{lem:inverse-shift-graph}
Let $\mathcal A,\mathcal B\subseteq\binom{[n]}k$ and suppose that
$(\mathcal U,\mathcal V)=(s_{ij}\mathcal A,s_{ij}\mathcal B)$
is a pair of cross $t$-intersecting families. On each side, take
the $(k-1)$-sets $P$ avoiding $i,j$ for which $P+i$ belongs to that
side and $P+j$ does not. Form a bipartite graph on the two collections,
joining $P$ and $Q$ if $|P\cap Q|=t-1$.
A pair with image $(\mathcal U,\mathcal V)$ is obtained by choosing
$P+i$ or $P+j$ at each vertex and fixing the other shift orbits.
Its families are cross $t$-intersecting if and only if the chosen
endpoint is constant on each connected component of the graph.
\end{lemma}
\begin{proof}
An inverse shift fixes an orbit containing both sets or neither set. A singly occupied orbit has the two stated choices. Between
two singly occupied orbits, the intersection decreases by one precisely
when their choices differ. Since the image families are cross
$t$-intersecting, the two bases meet in at least $t-1$ points.
Thus different choices violate cross $t$-intersection exactly
when the two bases are adjacent. A doubly occupied orbit
already contains both endpoint choices in
the image pair. Intersection with a set containing both endpoints,
or neither, is independent of the choice. Hence there are no further
restrictions.
\end{proof}

\begin{lemma}\label{lem:inverse-radius-one}
If $(s_{ij}\mathcal A,s_{ij}\mathcal B)=(\mathcal A_1(Z),\mathcal A_1(Z))$
and $\mathcal A,\mathcal B$ are cross $t$-intersecting, then
$\mathcal A=\mathcal B=\mathcal A_1(Z')$ for some $(t+2)$-set $Z'$.
\end{lemma}
\begin{proof}
If $i,j$ both belong to $Z$, or both avoid it, there are no singly
occupied orbits. If $j\in Z$ and $i\notin Z$, choose a $k$-set avoiding $i$
whose trace on $Z$ has size $t+1$ and contains $j$. Replacing
$j$ by $i$ removes it from $\mathcal A_1(Z)$, so this case
cannot be a shifted image. In the remaining case $i\in Z$, $j\notin Z$, a singleton
base has the form
\[
 (Z\setminus\{i,u\})\cup R,\qquad
 u\in Z\setminus\{i\},\quad
 R\in\binom{[n]\setminus(Z\cup\{j\})}{k-t-1}.
\]
Two bases are adjacent exactly when their omitted labels differ and
their outside parts are disjoint. Any two left vertices have a
common right neighbour: since $t+1\ge3$, we can choose a label
different from both omitted labels, and since $n-t-3\ge3(k-t-1)$,
we can choose an outside $(k-t-1)$-set disjoint from both outside
parts. For any single right vertex, the same argument gives
a left neighbour, so the graph is connected. Lemma~\ref{lem:inverse-shift-graph} permits only
the original pair or its simultaneous transposition of $i,j$.
\end{proof}

\begin{lemma}\label{lem:restricted-disjointness-connected}
Let $r\ge3$, $|V|\ge3r$, and let $P,Q\subseteq V$ each have
size $r$ or $r+1$. The bipartite disjointness graph with sides
\[
 \mathcal U=\{R\in\tbinom Vr:R\cap Q\ne\varnothing,\ P\nsubseteq R\},
 \qquad
 \mathcal V=\{S\in\tbinom Vr:S\cap P\ne\varnothing,\ Q\nsubseteq S\}
\]
is connected.
\end{lemma}
\begin{proof}
For each $p\in P$, consider the left vertices avoiding $p$.
Any two such vertices $R,R'$ have a common right neighbour: choose
an $r$-set in $V\setminus(R\cup R')$ containing $p$.
There are at least $r$ available points, and the chosen set cannot
contain all of $Q$, since $R$ meets $Q$.
For distinct $p,p'\in P$, the two classes of left vertices intersect:
choose an $r$-set avoiding $p,p'$ and containing a point of
$Q\setminus\{p,p'\}$, which is nonempty as $|Q|\ge3$.
Every left vertex avoids some point of $P$, so all left vertices lie
in one component. Finally, every right vertex $S$ misses some
$q\in Q$. Choose $p\in P\setminus\{q\}$, and extend $\{q\}$
to an $r$-set
in $V\setminus(S\cup\{p\})$. This is possible since
$|V\setminus S|\ge2r$. The resulting set is a left neighbour of $S$.
\end{proof}

\begin{lemma}\label{lem:inverse-dual-hm}
Let $r=k-t\ge3$, $|T|=t$, $|X|=|Y|=r+1$,
$T\cap(X\cup Y)=\varnothing$ and $|X\cap Y|\ge2$.
Suppose
\[
 (s_{ij}\mathcal A,s_{ij}\mathcal B)
   =(\mathcal H(T;X,Y),\mathcal H(T;Y,X))
\]
and $\mathcal A,\mathcal B$ are cross $t$-intersecting.
Then $(\mathcal A,\mathcal B)$ is also a dual Hilton--Milner pair
whose two outside blocks have size $r+1$ and meet in at least two
points.
\end{lemma}
\begin{proof}
Put $c=|X\cap Y|$. We analyse the graph of
Lemma~\ref{lem:inverse-shift-graph}.
If $i,j\in T$, both families are invariant under their transposition.
The case $j\in T$, $i\notin T$ cannot occur. Fix
$y\in Y\setminus\{i\}$ and choose an $r$-set $R$ outside
$T\cup\{i\}$ that contains $y$ and satisfies $R\cup\{i\}\ne X$.
There are at least two possible choices for $R$ containing $y$,
and at most one is excluded. The member $T\cup R$ would shift
to $(T\setminus\{j\})\cup(R\cup\{i\})$, which is absent
from the first family. Thus the displayed image would not be
$ij$-compressed.

First suppose $i,j\notin T$. Since the image families are
$ij$-compressed,
\begin{equation}\label{eq:inverse-outside-memberships}
 j\in X\Longrightarrow i\in X,\qquad
 j\in Y\Longrightarrow i\in Y.
\end{equation}
For example, if $j\in X$ and $i\notin X$, shifting any exceptional
member $(T\setminus\{p\})\cup X$ changes its outside block, and
the resulting set is absent from that family. If neither block
contains $i$ without $j$, there are no singleton orbits.

Suppose that $i\in X$, $j\notin X$, whereas $Y$ contains both
endpoints or neither. The singleton bases on
the first side are
$(T\setminus\{p\})\cup(X\setminus\{i\})$, $p\in T$;
those on the second are the main-star bases $T\cup R$ with
$R\in\binom{[n]\setminus(T\cup X\cup\{j\})}{r-1}$.
Both collections are nonempty, since $n>2k$. Every base in the first
collection is adjacent to every base in the second. The common choice of endpoint either retains $X$ or replaces
$i$ by $j$ in $X$. By~\eqref{eq:inverse-outside-memberships}, $Y$
contains both endpoints or neither, so the size of the intersection
of the two blocks is unchanged.
The case where only $Y$ contains $i$ without $j$ is symmetric.

If $i\in X\cap Y$ and $j\notin X\cup Y$, the singleton bases
form two complete bipartite subgraphs: the main-star bases on the first
side with the exceptional bases on the second, and conversely.
Two exceptional bases with different omitted points of $T$ meet
in $t+c-3$ points. Since $t\ge2$, there is an edge between the
two subgraphs when $c=2$.
When $c\ge3$, they have no edges between them: two main-star bases
already share $T$, and two exceptional bases intersect in at least
$t$ points.
In the first case both blocks must move together; in the second each
block can independently retain $i$ or replace it by $j$.
Moving just one block decreases the size of their intersection
from $c$ to $c-1$, while moving both preserves it. The two blocks therefore still meet in at least two points.
A change of the endpoint in an exceptional block acts on exactly
the corresponding main-star singleton orbits of the opposite
family, namely those whose outside part otherwise misses the
block. Thus the inverse pair has the claimed form with the two
new outside blocks.

It remains that $i\in T$ and $j\notin T$. Put
$T_0=T\setminus\{i\}$, $V=[n]\setminus(T\cup\{j\})$,
$X_0=X\setminus\{j\}$ and $Y_0=Y\setminus\{j\}$.
The main-star singleton bases on the first side are $T_0\cup R$,
where
\[
 R\in\binom Vr,\qquad R\cap Y_0\ne\varnothing,
 \qquad X_0\nsubseteq R.
\]
If $j\in X$, the exceptional member $T_0\cup X$ and the main-star
member $T\cup(X\setminus\{j\})$ occupy both positions in one orbit:
the latter is a member since $|X\cap Y|\ge2$. The last restriction
excludes this orbit; if $j\notin X$, it is automatic.
The analogous second-side bases interchange $X_0,Y_0$.
Their graph is precisely that of
Lemma~\ref{lem:restricted-disjointness-connected}, since their common
$T_0$ has size $t-1$ and $|V|=n-t-1\ge(t+1)r\ge3r$.
The block sizes are $r$ or $r+1$, so the lemma gives connectivity.

There are additional exceptional singleton bases only if $j\notin X$
on the first side or $j\notin Y$ on the second. A first-side one is
$(T_0\setminus\{p\})\cup X$, $p\in T_0$.
To find a second-side main-star base $T_0\cup R$ with $|R\cap X|=1$,
take $R$ to consist of one point of $X$ and $r-1$ points of
$V\setminus X$. This is possible because $|V|\ge3r$ and $|X|=r+1$.
It cannot equal $Y_0$ when $|Y_0|=r$, since
$|Y_0\cap X|=|Y\cap X|\ge2$ when $j\notin X$.
The two bases have intersection size $t-1$, so this exceptional
vertex has a neighbour in the connected main-star graph. The same
argument applies to exceptional vertices on the other side. Thus
the entire graph is connected, and the inverse pair is either the
image pair or its simultaneous transposition of $i,j$.
\end{proof}

\subsection{Proof of the main theorem}\label{sec:main-proof}
\begin{proof}[Proof of Theorem~\ref{thm:main}]
Among all pairs of nontrivial cross $t$-intersecting families with
maximum product, choose a pair that minimizes the sum of the element
labels over both families. This pair is saturated. Any simultaneous
left shift that changes the pair must make at least one family
trivial. If there is such a shift, Theorem~\ref{thm:critical} gives
a product strictly less than $m_t(n,k)^2$. Otherwise both families
are shifted, and Proposition~\ref{prop:shifted-bound} gives the
required bound. In either case the bound holds for a pair of
maximum product and hence for every admissible pair.
The Hilton--Milner and radius-one constructions have sizes
$h_t(n,k)$ and $b_t(n,k)$, respectively. Pairing a construction
of maximum size with itself attains the bound.

Now let $\mathcal A,\mathcal B$ be nontrivial cross $t$-intersecting
families with $|\mathcal A||\mathcal B|=m_t(n,k)^2$. Apply simultaneous
shifts that keep both families nontrivial until the pair is terminal.
Every step preserves both sizes. By Theorem~\ref{thm:critical},
the terminal pair is shifted, since otherwise its product would
be smaller. Lemmas~\ref{lem:shifted-equality-cover}
and~\ref{lem:transversal-equality} show that the shifted pair is
of type \textup{(i)} or \textup{(ii)}. Reverse the shifts one at a time. The radius-one class is
preserved by Lemma~\ref{lem:inverse-radius-one}, and the dual
Hilton--Milner class by Lemma~\ref{lem:inverse-dual-hm}.
The latter lemma requires $k-t\ge3$, which is automatic whenever
$h_t(n,k)\ge b_t(n,k)$: for $k=t+2$,
$b_t(n,k)-h_t(n,k)=(t-1)(n-t-3)>0$.
This proves necessity.

Conversely, the radius-one pair in \textup{(i)} is nontrivial cross
$t$-intersecting and has the required product. In \textup{(ii)},
each family has size $h_t(n,k)$. Two members from opposite main
stars share $T$. A main-star member and an opposite exceptional
member share $t-1$ points of $T$ and at least one outside point.
Two exceptional members share at least $t-2+|X\cap Y|\ge t$ points.
The intersection of the main-star members is $T$, and for each
point of $T$ an exceptional member omits that point. Thus both
families are nontrivial, and the pair attains the bound.
\end{proof}

\appendix
\section{Algebraic inequalities for generator exchanges}\label{sec:exchange-numerics}

We first record the calculation common to the four lemmas below.
For $t\ge2$, $t+1\le i\le j\le k$, $s=i+j-t$ and $n>2k-t+1$, put
$C_v=\binom{n-s}{k-v}$ and $C_v^+=\binom{n-s}{k-v+1}$.
Consider the paired exchange products
\begin{equation}\label{eq:appendix-products}
 (A+\alpha C_i^+)(B-\beta C_j),\qquad
 (A-\alpha C_i)(B+\beta C_j^+),
 \qquad A,B,\alpha,\beta>0.
\end{equation}
Write $U=A/(\alpha C_i)$, $V=B/(\beta C_j)$ and $R_v=C_v^+/C_v$.
For $x=k-i+1$, $y=k-j+1$ and $h=n-2k+t-1>0$, we have
$R_i=(h+y)/x$ and $R_j=(h+x)/y$, so $R_iR_j>1$.
If both products are at most $AB$, then
$R_iV\le U+R_i$ and $R_jU\le V+R_j$. Eliminating $V$ gives
\begin{equation}\label{eq:appendix-upper}
 U\le\frac{R_i(R_j+1)}{R_iR_j-1}
   =1+\frac{k-j+1}{n-2k+t-1}.
\end{equation}
The shadow bound
\begin{equation}\label{eq:appendix-shadow}
 U\ge1+\frac{s-i}{i}\frac{n-s+1}{n-k-s+i+1}
\end{equation}
contradicts this upper bound whenever
\begin{equation}\label{eq:appendix-polynomial}
 P(n):=(j-t)(n-i-j+t+1)(n-2k+t-1)
       -i(k-j+1)(n-k-j+t+1)>0.
\end{equation}
For balanced ranks, \eqref{eq:appendix-upper} reads $U\le R_i/(R_i-1)$.

\begin{lemma}\label{lem:general-individual-exchange}
Let $t\ge3$, $s\ge t+3$ and $n\ge(t+1)(k-t+1)$.
Under the notation above and the shadow bound~\eqref{eq:appendix-shadow},
at least one product in~\eqref{eq:appendix-products} exceeds $AB$.
\end{lemma}
\begin{proof}
Write $i=t+p$, $j=t+q$, $k=t+r$, $n_0=(t+1)(r+1)$ and $z=r-q$.
Then $1\le p\le q\le r$, $q\ge2$ and $z\ge0$.
Since $P''(n)=2q>0$ and $P'(n_0)$ decreases with $p$,
\[
 P'(n_0)\ge(2qt-q-t)z+2q^2(t-1)+q-t>0.
\]
It remains to check $P(n_0)>0$. This too decreases with $p$, so
\begin{align*}
 P(n_0)\ge{}&[q(t^2-t-1)-t^2]z^2\\
 &+[(2t^2-3t+1)(q-2)^2+(7t^2-10t)(q-2)+5t^2-10t-4]z\\
 &+(q+1)(qt-q+2)(qt-q-t)>0.
\end{align*}
All coefficients are positive for $t\ge3$, $q\ge2$.
Thus \eqref{eq:appendix-polynomial} holds for $n\ge n_0$.
\end{proof}

\begin{lemma}\label{lem:t2-individual-exchange}
Let $t=2$, $q=s=i+j-2$ and $n\ge3k-3$.
Under~\eqref{eq:appendix-shadow}, at least one product
in~\eqref{eq:appendix-products} exceeds $AB$ if $i<j$ or $i=j\ge6$.
\end{lemma}
\begin{proof}
Put $z=k-j\ge0$. We have $P''(n)=2(j-2)>0$ and
$P'(3k-3)=(4j-8-i)z+(j-2)(3j-i-2)-i>0$. Moreover,
\[
 P(3k-3)=(3j-6-2i)z^2+(5j^2-2ij-16j+12)z+(j-2)^2(2j-i)-ij.
\]
For $j\ge i+1$, the coefficients increase with $j$ and at $j=i+1$
equal $i-3$, $3i^2-8i+1$ and $i^3-i^2-4i+2$.
For $j=i\ge6$, they equal $i-6$, $3i^2-16i+12$ and $i(i-1)(i-4)$.
In either case the first is nonnegative and the others are positive,
proving~\eqref{eq:appendix-polynomial}.
\end{proof}

The remaining balanced ranks are $i=j=3,4,5$.
The rank-three case is a terminal configuration in the main proof;
the next two lemmas treat ranks four and five.

\begin{lemma}\label{lem:t2-global-middle}
Let $(q,i)\in\{(6,4),(8,5)\}$, $n\ge3k-3$ and $k\ge i$.
Let $\mathcal A,\mathcal B$ be nontrivial shifted cross $2$-intersecting
families whose canonical generators are supported on $[q]$.
Suppose the paired rank-$i$ exchange has nonempty active layers and
the increments in Lemma~\ref{lem:paired-terminal-exchange}.
If either family has a generator of size at most $i-1$, at least one
exchange strictly increases the product.
\end{lemma}
\begin{proof}
By symmetry, assume $\mathcal A$ has a generator of size at most $i-1$.
If this generator contains $q$, shift $q$ to a missing earlier point
and choose a canonical generator inside the resulting set.
This gives $F\subseteq[q-1]$ with $h=|F|\le i-1$.
The extension argument in Lemma~\ref{lem:canonical-pairing}
shows that $F$ meets every canonical generator of $\mathcal B$ in at
least two points. Since $\mathcal B$ is nontrivial, $h\ge3$.
Let $\mathcal T$ be the traces on $[q-1]$ of its active rank-$i$
generators, put $\alpha=|\mathcal T|$, and write $\nabla\mathcal T$
for their rank-$i$ upper shadow. By minimality of the generators,
no trace contains $F$. Since $h\le i-1$, at least one $(i-1)$-set
containing $F$ is missing from $\mathcal T$.
For $(q,i)=(6,4)$ set $(c,w)=(5/9,1/9)$; for $(q,i)=(8,5)$ set
$(c,w)=(21/34,1/34)$. Then
\begin{equation}\label{eq:proper-small-shadow}
 \alpha\le1/w,\qquad |\nabla\mathcal T|\ge c\alpha.
\end{equation}
For the shadow estimate, take complements and apply
Kruskal--Katona~\cite{Kruskal,Katona}. With at most nine edges on five
points, colex has shadow-to-size ratio at least $2/3$ up to six edges,
and uses five vertices thereafter. With at most 34 triples on seven
points, the ratio is at least $3/4$ up to 20 triples. For 21--30
triples it is at least $(15+v)/(20+\binom v2)\ge2/3$, where
$2\le v\le5$; for 31--34 triples the shadow contains all 21 pairs.

Put $L=n-q$, $C=\binom L{k-i}$, $R=\binom L{k-i+1}/C$ and
$H=\binom{L+1}{k-i}/C$. Besides the active and shadow cylinders,
$\mathcal A$ contains the cylinder with trace $F$ on $[q-1]$.
It is disjoint from the former because $F\notin\mathcal T$, and
from the latter because $h<i$. Its size is
$\binom{L+1}{k-h}\ge\binom{L+1}{k-i+1}=(R+1)C$:
indeed, $k-i+1\le k-h$ and $2(k-h)\le2k-6\le L+1$. Consequently
\[
 \frac{|\mathcal A|}{\alpha C}\ge1+cH+w(R+1).
\]
This exceeds the obstruction $R/(R-1)$ from~\eqref{eq:appendix-upper}.
Indeed, with $D_1=n-k-q+i+1>D_2=n-2k+1>0$, the difference has derivative
\[
 -\frac{c(k-i)}{D_1^2}+\frac{w}{k-i+1}
       +\frac{k-i+1}{D_2^2}>0.
\]
At $n=3k-3$, its values in the two cases are
\[
 \frac{3k^2-5k-10}{18(k-3)(k-2)}>0,\qquad
 \frac{(k-5)^3+94(k-5)^2+301(k-5)+220}
      {34(k-4)(k-2)(2k-5)}>0,
\]
respectively. Thus one exchange increases the product.
\end{proof}

Under the hypotheses of Lemma~\ref{lem:t2-global-middle}, if the pair
is saturated and maximizes the product and both exchanges preserve
nontriviality, then every generator has size at least $i$. Every member
on either side therefore has at least $i$ points in $[q]$. Since
$q=2i-2$, all $k$-sets with this property meet every member of either
family in at least two points, so saturation gives
$\mathcal A=\mathcal B=\{A\in\binom{[n]}k:|A\cap[q]|\ge i\}$.
These radius-two and radius-three families have sizes $F_2$ and $F_3$,
respectively, in the proof of Lemma~\ref{lem:asymmetric-threshold},
where $b_2(n,k)=F_1>F_2\ge F_3$.

\begin{lemma}\label{lem:t2-inner-extra-cylinder}
Let $(q,i)\in\{(6,4),(8,5)\}$, $k\ge q$ and $n\ge3k-3$.
Suppose a paired inner exchange has products~\eqref{eq:appendix-products}
with $t=2$, $s=q$ and $j=i$. Suppose the first family contains the
$\alpha$ active cylinders, their upper-shadow cylinders counted
in~\eqref{eq:terminal-shadow-bound}, and every $k$-set with trace
$[2]$ on $[q]$ that meets $[q+1,k+1]$.
Then at least one exchange product exceeds $AB$.
\end{lemma}
\begin{proof}
Put $L=n-q$, $R=C_i^+/C_i$ and
$U_0=1+\frac{q-i}{i}\frac{n-q+1}{n-k-q+i+1}$.
The extra cylinders have trace $[2]$, so they are disjoint from the
active and shadow cylinders, whose traces have at least $i-1\ge3$
points. Retain $d=\min\{3,k+1-q\}$ fixed outer points. Counting the
extra sets by their first point among the $d$ chosen points and using
$\alpha\le\binom{q-1}{i-1}$ gives
\[
 \frac A{\alpha C_i}\ge U_0+E_d,\qquad
 E_d=\frac1{\binom{q-1}{i-1}}
      \sum_{h=1}^d\frac{\binom{L-h}{k-3}}{\binom L{k-i}}.
\]
We show that this exceeds $R/(R-1)$, contradicting~\eqref{eq:appendix-upper}.

The difference $U_0-R/(R-1)$ increases with $n$: writing
$D_1=n-k-q+i+1>D_2=n-2k+1>0$, its derivative is
$-\frac{(q-i)(k-i)}{iD_1^2}+\frac{k-i+1}{D_2^2}>0$.
Each summand of $E_d$ also increases, being a positive constant times
\[
 \prod_{v=0}^{h-1}\left(1-\frac{k-i}{L-v}\right)
       (L-k+i-h)_{i-3}.
\]
It therefore suffices to take $n=3k-3$.

For $(q,i)=(6,4)$ and $k\ge8$, put $x=k-3$, so $L=3x$ and
$R/(R-1)-U_0=(x-1)/(4x+4)<1/4$.
Cancellation gives
\[
 \frac{\binom{3x-1}{x}}{\binom{3x}{x-1}}\ge\frac43,\quad
 \frac{\binom{3x-2}{x}}{\binom{3x}{x-1}}\ge\frac89,\quad
 \frac{\binom{3x-3}{x}}{\binom{3x}{x-1}}\ge\frac12.
\]
Hence $E_3\ge49/180>1/4$. For $k=6,7$, the exact differences
$U_0+E_d-R/(R-1)$ are $11/360$ and $21/220$, respectively.

For $(q,i)=(8,5)$,
$R/(R-1)-U_0=(k^2-17k+40)/((k-2)(10k-25))<1/10$.
At $k=8$ this is $-16/165$. For $k\ge9$, put $x=k-4$, so $L=3x+1$.
The first two cylinder ratios equal
\[
 \frac{4(2x+1)}{3x+1}>\frac83,\qquad
 \frac{4(4x^2-1)}{3x(3x+1)}\ge\frac32.
\]
Thus $E_d\ge(8/3+3/2)/35=5/42>1/10$, completing the proof.
\end{proof}

\end{document}